\documentclass[11pt]{article}
\usepackage{amsfonts}
\usepackage{amsfonts,latexsym,amsmath,amscd,geometry}
\usepackage{color}
\usepackage{amssymb}
\usepackage{latexsym}

\newcommand \nc{\newcommand}
\newtheorem{theorem}{Theorem}[section]
\newtheorem{lemma}[theorem]{Lemma}
\newtheorem{proposition}[theorem]{Proposition}
\newtheorem{corollary}[theorem]{Corollary}
\newtheorem{definition}[theorem]{Definition}

\newtheorem{remark}[theorem]{Remark}

\nc{\ba}{\begin{array}}\nc{\ea}{\end{array}}
\nc{\be}{\begin{eqnarray}}\nc{\ee}{\end{eqnarray}}
\nc{\beq}{\begin{equation}}\nc{\eeq}{\end{equation}}
\nc{\bex}{\begin{eqnarray*}}\nc{\eex}{\end{eqnarray*}}
\nc{\btm}{\begin{theorem}} \nc{\etm}{\end{theorem}}
\nc{\blm}{\begin{lemma}} \nc{\elm}{\end{lemma}}
\nc{\R}{\mathbb{R}} \nc{\va}{\varepsilon} \nc{\ls}{\limits}

\def\pf{\noindent{\bf Proof.\quad}}\def\endpf{\hfill$\Box$}

\begin{document}
\title{Global weak solution to
the viscous two-fluid model with finite energy}
\author{Alexis Vasseur\footnote{Department of Mathematics, The University
of Texas at Austin, USA; vasseur@math.utexas.edu.}\quad Huanyao
Wen \footnote{School of Mathematics, South China University of
Technology, Guangzhou, China; mahywen@scut.edu.cn.} \quad Cheng
Yu\footnote{Department of Mathematics, The University of Texas at
Austin, USA; yucheng@math.utexas.edu.}}

\maketitle

\begin{abstract}
In this paper, we  prove the existence of global weak solutions to
the compressible two-fluid Navier-Stokes equations in three dimensional space. The pressure depends on two different variables from the continuity equations.
 We develop an argument of variable reduction for the  pressure law. This yields to the strong convergence of the
densities, and provides the existence of global solutions in time,
for the compressible two-fluid Navier-Stokes equations, with
large data in three dimensional space.
\end{abstract}


{\bf \noindent keywords.} two-fluid model, compressible Navier-Stokes equations, global weak solutions.


\section{Introduction}
\setcounter{equation}{0} \setcounter{theorem}{0} In this paper, we are considering
a viscous compressible
two-fluid model with a pressure law in two
variables. We show  the existence of
global weak solutions to the following two-fluid compressible Navier-Stokes system in three dimensional space:
    \begin{equation}\label{equation}
    \left\{
    \begin{array}{l}
        n_t+\mathrm{div}(n u)=0,\\
        \rho_t+\mathrm{div}(\rho u)=0,\\
 \big[(\rho+n) u\big]_t+\mathrm{div}\big[(\rho+n) u\otimes u\big]
 +\nabla P(n,\rho)=
            \mu \Delta u + (\mu+\lambda)\nabla \mathrm{div}u \quad \mathrm{on} \   \Omega\times (0, \infty),
    \end{array}
    \right.
    \end{equation}
with  the  initial and boundary conditions

    \begin{equation}\label{initial}
      n(x,0)=n_0(x),\ \rho(x,0)=\rho_0(x),\ (\rho+n) u(x,0)=M_0(x)  \ \ \mathrm{for} \ \  x\in\overline{\Omega},
    \end{equation}

     \begin{equation}\label{boundary}
     u|_{\partial\Omega}=0\quad \mathrm{for}\ t\ge0,
    \end{equation}
    where $\Omega\subset\R^3$ is a bounded domain, $P(n,\rho)=\rho^{\gamma}+n^{\alpha}$ denotes the pressure for $\gamma\geq1\text{ and }\alpha\geq1$, $u$ stands for the velocity of fluid,
  $\rho$ and $n$ are the densities of two fluids, $\mu$ and $\lambda$ are the viscosity coefficients. Here we assume that $\mu\text{ and } \lambda$ are fixed constants, and $$\mu> 0,\;\;2\mu+\lambda\geq 0.$$

The two-fluid model was originally developed by Zuber and Findlay
\cite{Z}, Wallis \cite{W}, and Ishii \cite{I,I2}. The case
$\alpha=1$  corresponds to  the
hydrodynamic equations of \cite{CG, MV3}. It was  derived in  \cite{CG,
MV3} as an asymptotic limit of a
 coupled system of  the compressible Navier-Stokes equation with
a Vlasov-Fokker-Planck equation. The case $\alpha=2$ is associated to the  compressible Oldroyd-B type model with stress diffusion, see  Barrett, Lu, and Suli \cite{BLS}. The main difference with  the
classical compressible Navier-Stokes equations is that the pressure
law $P(\rho,n)=\rho^{\gamma}+n^{\alpha}$ depends on  two variables. In
this context, the existence of weak solutions to equations
(\ref{equation}) remained open until now.
 We refer the reader to
\cite{BBCMT,BCT, Br, I,I2,W,Z} for more physical background and
discussion of numerical studies for such mathematical models.

\vskip0.3cm

One difficulty dealing with the compressible Navier-Stokes equation is the degeneracy of the system close to the vacuum (when the density is vanishing).
The first existence result for the compressible Navier-Stokes equations in one dimensional space was established by Kazhikhov and Shelukhin \cite{KS}. This result was restricted to   initial densities bounded away from zero.
It has been extended by Hoff \cite{Hoff87}  and Serre \cite{S} to the case of   discontinuous initial data, and by Mellet-Vasseur \cite{MV2} in the case of density dependent viscosity coefficients.  For the multidimensional  case,
 the first global existence with  small initial data was proved by Matsumura and Nishida \cite{MN79,MN80,MN83}, and later by Hoff \cite{H95JDE,H95,H97} for discontinuous initial data.
 Lions, in \cite{Lions}, introduced the concept of  renormalized solutions for the compressible Navier-Stokes equations which allows to  control the possible oscillations of density. He proved
 the global existence  of 3D solutions for $\gamma\geq\frac{9}{5}$, and large initial values. It was  later improved
  by Jiang and Zhang \cite{Jiang-Zhang} for spherically symmetric initial data for $\gamma>1$, and by Feireisl-Novotn\'{y}-Petzeltov\'{a} \cite{Feireisl} and Feireisl \cite{Feireisl2}
  for $\gamma>\frac{3}{2}$,
  and to Navier-Stokes-Fourier systems. One key ingredient of the theory \cite{Lions, Jiang-Zhang, Feireisl} is to obtain
 higher
integrability on the density. This is obtained thanks to the elliptic structure on the viscous effective flux, and the specific form of the pressure  $P=\rho^{\gamma}.$ Relying on this structure, Lions
deduced that the density $\rho$ is uniformly bounded in
$L^{\gamma+\frac{2\gamma}{3}-1}$. Note that  for $1\leq \gamma\leq
\frac{3}{2}$, the construction of  weak solutions for large
data remains largely open, see \cite{Lions2}.  The primary
difficulty is the possible concentration of the convective term
in this case. Very recently, Hu \cite{Hu} studied the
concentration phenomenon of the kinetic energy, $\rho|u|^2$,
associated to the isentropic compressible Navier-Stokes equations for
$1\leq\gamma\leq\frac{3}{2}.$  Finally, let us mention a very promising work of
 Bresch-Jabin \cite{BJ}.  They developed  a new method to obtain compactness on the density. This method is very different from the theory initiated by Lions. It allows already the treatment of  non-monotone  pressure laws.
\vskip0.3cm

 The problem becomes even more challenging when the pressure law depends on two variables as follows \begin{equation}
 \label{nonlinear form}
 P(\rho,n)=\rho^{\gamma}+n^{\alpha}.
 \end{equation}
  To the best of our
knowledge, the only  results on global existence of
weak solutions to  System \eqref{equation} with large initial data are restricted to the one dimension case, see \cite{EK,EWZ} (See also
\cite{E,EK2,Yao-Zhang-Zhu} for smallness assumptions). In \cite{BLS},
Barrett-Lu-Suli established the existence of weak solutions to a
compressible Oldroyd-B type model with  pressure law
$P=\rho^{\gamma}+n+n^2$, in the two dimensional space, but with an extra diffusion term on the $n$ equation.  This provides higher regularity on $n$ due to
the parabolic structure.
David-Michalek-Mucha-Novotny-Pokorny-Zatorska  in \cite{MMMNPZ} constructed a weak solution of the compressible Navier-Stokes system with the nonlinear pressure law
$$P(\rho,s)=\rho^{\gamma}\mathcal{T}(s),\;\;\;\gamma\geq\frac{9}{5},$$ where $s$ satisfies the entropy equation.  Note that the quantity $\theta=(\mathcal{T}(s))^{\frac{1}{\gamma}}$ can be interpreted as a potential temperature, thus the pressure could take the form $P=(\rho\theta)^{\gamma}=Z^{\gamma}.$ The quantity $Z$ also satisfies the continuity equation. This allowed them    to apply the  standard technique for the compressible Navier-Stokes equations to this system. 

Because the pressure law  depends genuinely on two variables, the treatment of the system \eqref{equation}  is more involved.  At first sight, it seems that more regularity on the  densities is required  to control the cross products, like $\rho^{\gamma}n$ and $n^{\alpha}\rho$.  These extra regularity properties are, so far, out of reach, and the classical techniques cannot be applied directly on  \eqref{equation}.


\vskip0.3cm

For any smooth solution of system \eqref{equation},  the following energy inequality holds for any time $0\leq t\leq T:$
 \be\label{0-energy-inequality}
\begin{split} &\frac{d}{dt}\int_\Omega\Big[\frac{(\rho+n)|u|^2}{2}+G_\alpha(n)+\frac{1}{\gamma-1}\rho^\gamma\Big]\,dx
+\int_\Omega\Big[\mu|\nabla u|^2
+(\mu+\lambda)|\mathrm{div}u|^2\Big]\,dx\leq 0,
\end{split}
\ee  where \bex
G_\alpha(n)=\left\{\begin{array}{l} n\ln
n-n+1, \ \ \ {\rm for} \ \ \ \alpha=1, \\
[3mm] \frac{n^\alpha}{\alpha-1}, \ \ \ {\rm for}\ \ \ \alpha>1.
\end{array}
\right.\eex As usual, we assume that \begin{equation*}
\label{initial energy}
\int_\Omega\Big[\frac{(\rho_0+n_0)|u_0|^2}{2}+G_\alpha(n_0)+\frac{1}{\gamma-1}\rho_0^\gamma\Big]\,dx<\infty
\end{equation*}
in the whole paper. Thus, we set the following restriction on the initial data
\begin{equation}
\label{initial restriction-1}
\begin{split}&
\inf\limits_{x\in \Omega}\rho_0\geq0,\quad \inf\limits_{x\in
\Omega}n_0\geq0,\quad \rho_0\in L^\gamma(\Omega),\quad
G_\alpha(n_0)\in
 L^1(\Omega),
 \end{split}
 \end{equation}
and
\begin{equation}
\label{restriction-2}
 \frac{M_0}{\sqrt{\rho_0+n_0}}\in L^2(\Omega) \;\text{ where }\,\frac{M_0}{\sqrt{\rho_0+n_0}}=0\;\text{ on }\,\{x\in\Omega|\rho_0(x)+n_0(x)=0\}.
\end{equation}
The definition of weak solution in the energy space is given in the following sense.
\begin{definition}
\label{definition of weak soluton}
We call $(\rho, n, u):\Omega\times (0,\infty)\to\mathbb
R_+\times\mathbb R_+\times\mathbb R^3$ a global weak solution of
(\ref{equation})-(\ref{boundary}) if for any $0<T<+\infty$,
\begin{itemize}
\item $ \rho\in L^\infty\big(0,T;L^\gamma(\Omega)\big),\ G_{\alpha}( n)\in
L^\infty\big(0,T;L^1(\Omega)\big),\ \sqrt{\rho+n} u\in
L^\infty\big(0,T;L^2(\Omega)\big), u\in
L^2\big(0,T;H_0^1(\Omega)\big),$

\item  $(\rho, n, u)\ \mathrm{solves}\ \mathrm{the}\
\mathrm{system}\ (\ref{equation})\ \mathrm{in}\
\mathcal{D}^\prime(Q_T),\ \mathrm{where}\
\mathrm{Q_T}=\Omega\times(0,T),$

\item $\big(\rho, n, (\rho+n) u\big)(x,0)=\big(\rho_0(x),n_0(x),
M_0(x)\big),\quad \mathrm{for\ a.e.}\ x\in \Omega,$

\item The energy inequality \eqref{0-energy-inequality} holds in
  $\mathcal{D}^\prime\big(\mathbb
R^3\times (0,T)\big), $

\item $(\ref{equation})_1\ \mathrm{and} (\ref{equation})_2\
\mathrm{hold}\ \mathrm{in}\ \mathcal{D}^\prime\big(\mathbb
R^3\times (0,T)\big)\ \mathrm{provided}\ \rho,n,u\ \mathrm{are}\
\mathrm{prolonged}\ \mathrm{to}\ \mathrm{be}\ \mathrm{zero}\
\mathrm{on}\ \mathbb R^3/\Omega,$

\item  the equation (\ref{equation})$_1$ and (\ref{equation})$_2$
are satisfied in the sense of renormalized solutions, i.e.,\bex
\partial_tb(f)+{\rm{div}}\big(b(f)u\big)+[b^\prime(f)f-b(f)]{\rm{div}}u=0\eex
holds in  $\mathcal{D}^\prime(Q_T)$, for any $b\in C^1(\mathbb R)$
such that $b^\prime(z)\equiv0$ for all $z\in \mathbb R$ large
enough, where $f=\rho,n$.
\end{itemize}
\end{definition}

The  main result of this paper is as follows.
\begin{theorem} \label{th:1.1} Assume that $\Omega\subset \R^3$ is a bounded domain in $\R^3$ of class $C^{2+\nu},$ $\nu>0$.
Let the initial data be under the conditions \eqref{initial
restriction-1}-\eqref{restriction-2}.
\begin{itemize}
\item{
 If  \bex
\alpha\geq1,\quad \gamma>\frac{9}{5}, \eex and the initial data additionally satisfies
\begin{equation}
\label{additional initial data-theorem}
 {\rm \ \frac{1}{c_0}\rho_0\le n_0\le c_0\rho_0\ on\
\Omega},
    \end{equation} where $c_0\ge1$ is a
constant,
 then there exists a global weak solution $(\rho, n, u)$ to
(\ref{equation})-(\ref{boundary}).}
\item{
Without restriction \eqref{additional initial data-theorem}, if \be\label{alphagamma}\alpha,\gamma>\frac{9}{5}\,\,\,
\mathrm{and}\,\,\, \max\{\frac{3\gamma}{4},\gamma-1,\frac{3(\gamma+1)}{5}\}<\alpha<\min\{\frac{4\gamma}{3},\gamma+1,\frac{5\gamma}{3}-1\},
\ee
 then there exists a global weak solution $(\rho, n, u)$ to
(\ref{equation})-(\ref{boundary}).}

\end{itemize}
\end{theorem}

\vskip0.3cm
\begin{remark}
The restriction  $\gamma>\frac{9}{5}$ provides the $L^2-$estimate on the density. This is needed to  apply the renormalized argument of DiPerna-Lions  to the system. The condition  \eqref{additional initial data-theorem} is propagated in time, and gives that
$$n\leq C\rho$$ for almost every time $t>0$. This provides extra integrability on $n$ without more assumption on the system  than  $\alpha\geq 1.$ However, without the condition \eqref{additional initial data-theorem}, the value of $\alpha$ needs to be close enough to the value of $\gamma$. This is required   to insure the  $L^2-$estimate on $n$. \end{remark}

The key idea of our proof is to perform a  variable reduction in the pressure law.  When considering a family of solutions, we decompose the pressure as
$$P_{\varepsilon}=\rho_{\varepsilon}^{\gamma}+n_{\varepsilon}^{\alpha}=A^\alpha d_{\varepsilon}^{\gamma}+B^\gamma d_{\varepsilon}^{\alpha}+\text{remainder},$$
where $d_{\varepsilon}=\rho_{\varepsilon}+n_{\varepsilon}.$
The idea is that we can control the oscillations of
$A_{\varepsilon}=n_{\varepsilon}/d_{\varepsilon}$ and $B_{\varepsilon}=\rho_{\varepsilon}/d_{\varepsilon}$.

\vskip0.3cm The structure of this paper is  as follows.
 In Section 2, we develop a new tool to handle the compactness on the terms $A_\varepsilon$ and $B_\varepsilon$.
 In section 3, we solve the approximation system using  the  Galerkin method.
In section 4, we study the limits as $\epsilon$ goes to zero. The
focus of this section is to  prove that
$$\overline{n^\alpha+\rho^\gamma+\delta(\rho+n)^\beta}=n^\alpha+\rho^\gamma+\delta(\rho+n)^\beta.$$
Here and always, $\overline{f}$ is the weak limit of
$f_{\epsilon}$. One of the key step is to control the product of $
n_\epsilon^\alpha+\rho_\epsilon^\gamma$ and $
n_\epsilon+\rho_\epsilon$, we rewrite them as follows \bex
\begin{split}
n_\epsilon^\alpha+\rho_\epsilon^\gamma=&A_\epsilon^\alpha
d_\epsilon^\alpha+B_\epsilon^\gamma d_\epsilon^\gamma =A^\alpha
d_\epsilon^\alpha+B^\gamma
d_\epsilon^\gamma+(A_\epsilon^\alpha-A^\alpha)
d_\epsilon^\alpha+(B_\epsilon^\gamma-B^\gamma) d_\epsilon^\gamma,\\
n_\epsilon+\rho_\epsilon=&(A_\epsilon+B_\epsilon) d_\epsilon
=(A+B)d_\epsilon +(A_\epsilon-A+B_\epsilon-B) d_\epsilon,
\end{split}
\eex where $d_\epsilon=\rho_\epsilon+n_\epsilon$, $d=\rho+n$,
$(A_\epsilon,B_\epsilon)=(\frac{n_\epsilon}{d_\epsilon},\frac{\rho_\epsilon}{d_\epsilon})$
if $d_\epsilon\neq0$, $(A,B)=(\frac{n}{d},\frac{\rho}{d})$ if
$d\neq0$, and $0\le A_\epsilon,B_\epsilon, A, B\le 1$, and
$(A_\epsilon d_\epsilon,B_\epsilon
d_\epsilon)=(n_\epsilon,\rho_\epsilon)$, $(Ad,Bd)=(n,\rho)$,
$(\rho,n)$ is the limit of $(\rho_{\epsilon},n_{\epsilon})$ in a
suitable weak topology.   Here we want to show
$$\Big[(A_\epsilon^\alpha-A^\alpha)
d_\epsilon^\alpha+(B_\epsilon^\gamma-B^\gamma)
d_\epsilon^\gamma\Big](n_\epsilon+\rho_\epsilon)\to 0,\,\,
\mathrm{and}\,\, \Big(A^\alpha d_\epsilon^\alpha+B^\gamma
d_\epsilon^\gamma\Big)(A_\epsilon-A+B_\epsilon-B) d_\epsilon\to
0$$ in some sense as $\epsilon\to 0.$ This can be done because
$\rho_\epsilon$ and $n_\epsilon$ are bounded uniformly for
$\epsilon$ in $L^{\beta+1}(Q_T)$ where
$\beta>\max\{\alpha,\gamma,4\}$, and, thanks to the result of Section 2
\begin{equation*}\label{key-equality}
\lim\limits_{\epsilon\rightarrow0^+}\int_0^T\int_\Omega
d_\epsilon|A_\epsilon-A|^s\,dx\,dt=0,\,\,\, \mathrm{and}\,\,\,
\lim\limits_{\epsilon\rightarrow0^+}\int_0^T\int_\Omega
d_\epsilon|B_\epsilon-B|^s\,dx\,dt=0.
\end{equation*}
In section 5, we recover the weak solution by letting $\delta$ goes to zero. The highlight of this section is to show that
 $$\overline{n^\alpha+\rho^\gamma}=n^\alpha+\rho^\gamma,$$
  which is similar to the limits in term of $\epsilon.$ However, a new difficulty occurs because of  $\delta=0$.  We follow \cite{Feireisl} and use
 a cut-off function in the renormalization to  show the strong convergence of $\rho_{\delta}$ and $n_{\delta}$. This can be done using again the variable reduction of Section 2.
At this level of approximation, we require
$\gamma>\frac{9}{5}$ such that $\rho_\delta$ is bounded in
$L^{\gamma+\theta_2}(Q_T)$ with $\gamma+\theta_2>2$ for some
$\theta_2$ satisfying $\theta_2<\frac{\gamma}{3}$ and
$\theta_2\le\min\big\{1,\frac{2\gamma}{3}-1\big\}$. In order to
guarantee that $n_\delta$ is bounded in $L^{q_1}(Q_T)$ for some
$q_1>2$, we require either
$\alpha\in(\frac{9}{5},\infty)\cap(\gamma-\theta_1,\gamma+\theta_2)$
or $\alpha\in[1,\infty)$ and $\frac{1}{c_0}\rho_0\le n_0\le
c_0\rho_0$.

\section{An error estimate} Our main goal of this section is to prove the following Theorem \ref{main 2-le:important}. The proof  relies on the
 the DiPerna-Lions theory of the renormalized
solutions to the transport equation.  This theorem allows us to obtain the weak stability of solutions to \eqref{equation}.
We start from the following lemma in this section.
\begin{lemma}\label{2-le:important}Let $\{(g_K,h_K)\}_{K=1}^{\infty}$ be a sequence with the following properties
\begin{equation}
\label{pro11}(g_K,h_K)\rightarrow (g,h)\ \mathrm{weakly}\
\mathrm{in}\ L^{p}(Q_T)\ \mathrm{as}\ K\to\infty,
\end{equation} for any given $p > 1$, $g_K,h_K\ge0$, and
\begin{equation}
\label{pro22} \lim\limits_{K\to+\infty}\int_0^T\int_\Omega
a_Kh_K\,dx\,dt\le\int_0^T\int_\Omega h a_h\,dx\,dt,
\end{equation}
where $a_K=\frac{h_K}{g_K}$ if
$g_K\neq0$, $a_h=\frac{h}{g}$ if $g\neq0$, $0\le a_K,a_h\le
\mathcal{C}$ for some positive constant $\mathcal{C}$ independent
of $K$, and $a_Kg_K=h_K$, $a_hg=h$, then
 \be\label{lelim00}
\lim\limits_{K\to+\infty}\int_0^T\int_\Omega
g_K|a_K-a_h|^2\,dx\,dt=0. \ee
 In particular,
\be\label{2-ik00} \lim\limits_{K\to+\infty}\int_0^T\int_\Omega
g_K|a_K-a_h|^s\,dx\,dt=0, \ee for any $s>1$.
\end{lemma}
\pf Note that \bex \int_0^T\int_\Omega
g_K|a_K-a_h|^2\,dx\,dt=\int_0^T\int_\Omega a_Kh_K\,dx\,dt
-2\int_0^T\int_\Omega h_K a_h\,dx\,dt+\int_0^T\int_\Omega g_K
a_h^2\,dx\,dt, \eex one obtains
 \bex\begin{split}
\lim\limits_{K\to+\infty}\int_0^T\int_\Omega
g_K|a_K-a_h|^2\,dx\,dt\le& \int_0^T\int_\Omega h a_h\,dx\,dt
-2\int_0^T\int_\Omega h a_h\,dx\,dt+\int_0^T\int_\Omega h
a_h\,dx\,dt\\=&0,
\end{split}
\eex where we have used $a_Kg_K=h_K$, $a_hg=h$, \eqref{pro22} and
the weak compactness of $g_K$ and $h_K$ in \eqref{pro11}. This
deduces (\ref{lelim00}).

By the H\"older inequality and (\ref{lelim00}), \eqref{2-ik00}
follows for  $s\in(1,2)$. If $s\in[2,\infty)$, note that
$(a_K-a_h)$ is bounded in
$L^\infty\big(\Omega\times(0,T)\big)$. This allows us to have
\eqref{2-ik00}.
\endpf\\
\newline
The following theorem is our main result of this section.

\begin{theorem}\label{main 2-le:important}
Let $\nu_K\to 0$ as $K\to+\infty,$ and $\nu_K\geq 0$. If
$\rho_K\geq 0$ and $n_K\geq 0$ are the solutions to
\begin{equation}
\label{the parabolic for density}
(\rho_K)_t+\mathrm{div}(\rho_Ku_K)=\nu_K\Delta\rho_K,\;\;\rho_K|_{t=0}=\rho_0,
\;\;\nu_K \frac{\partial\rho_K}{\partial \nu}|_{\partial\Omega}=0,
\end{equation}
and
\begin{equation}
\label{second parabolic equation for n}
(n_K)_t+\mathrm{div}(n_Ku_K)=\nu_K\Delta n_K,\;\;n_K|_{t=0}=n_0, \;\;\nu_K \frac{\partial n_K}{\partial \nu}|_{\partial\Omega}=0,
\end{equation}
respectively, with $C_0\ge1$ independent of $K$ such that
\begin{itemize}

\item   $
  \|(\rho_K,n_K)\|_{L^{\infty}(0,T;L^2(\Omega))}\leq C_0,\;\sqrt{\nu_K}\|\nabla\rho_K\|_{L^2(0,T;L^2(\Omega))}\leq C_0,\,\sqrt{\nu_K}\|\nabla n_K\|_{L^2(0,T;L^2(\Omega))}
  \leq C_0.
$
  \item $
  \|u_K\|_{L^2(0,T;H^1_0(\Omega))}\leq C_0.
$
\item for any $K>0$ and any $t> 0$:
\begin{equation}
\label{AAAinitial condition for n2 over
density}\int_{\Omega}\frac{b_K^2}{d_K}\,dx\leq
\int_{\Omega}\frac{b_0^2}{d_0}\,dx,
\end{equation}
\end{itemize}
where $b_K=\rho_K$ or $n_K$, and $d_K=\rho_K+n_K$.

Then, up to a subsequence, we have
$$n_K\to n,\;\;\;\rho_K\to\rho\;\;\text{ weakly in } L^\infty(0,T;L^2(\Omega)),$$
$$ u_K\to u \text{ weakly in } L^2(0,T;H^1_0(\Omega)),$$
and for any $s>1$,
\begin{equation}
\label{key bound}
\lim_{K\to+\infty}\int_0^T\int_{\Omega}d_K|a_K-a|^s\,dx\,dt=0,
\end{equation}
where $a_K=\frac{b_K}{d_K}$ if $d_K\neq0$, $a=\frac{b}{d}$ if
$d\neq0$, and $a_Kd_K=b_K$, $ad=b$. Here $(b,d)$ is the weak
limit of $(b_K,d_K)$.
\end{theorem}
\begin{remark}  The proof relies on the Diperna-Lions renormalized argument for  transport equation. The $L^2$ bounds of the densities $\rho_K$ and $n_K$ make it possible to use this theory  for equations \eqref{the parabolic for density}.
\end{remark}

\begin{remark}
\label{remark for an inequality} If $\nu_K>0$, $u_K$ is smooth
enough and $\rho_K$ is bounded by below, then \eqref{AAAinitial
condition for n2 over density} is verified. In fact, choosing
$\varphi(b_K,d_K)=\frac{b_K^2}{d_K},$ one obtains
\begin{equation*}
\begin{split}&
\frac{\partial \varphi(b_K,d_K)}{\partial t}+\mathrm{div}(\varphi
u_K ) +[\frac{\partial \varphi}{\partial b_K}b_K+\frac{\partial
\varphi}{\partial d_K}d_K-\varphi]\mathrm{div}u_K
\\&+\nu_K(\frac{\partial^2 \varphi}{\partial b_K^2}|\nabla b_K|^2
+\frac{\partial ^2\varphi}{\partial d_K^2}|\nabla
d_K|^2+2\frac{\partial ^2\varphi}{\partial b_K\partial d_K}\nabla
b_K\cdot\nabla d_K)-\nu_K\Delta \varphi=0.
\end{split}
\end{equation*}
Note that $$\frac{\partial \varphi}{\partial
b_K}b_K+\frac{\partial \varphi}{\partial d_K}d_K-\varphi=0$$ and
$\varphi$ is convex, thus we have
$$\frac{d}{dt}\int_{\Omega}\varphi(b_K,d_K)\,dx\leq 0.$$
\end{remark}
We will rely on the following lemma to show  Theorem \ref{main 2-le:important}.
\begin{lemma}
\label{main lemma}Let $\beta:\R^N\to \R$ be a $C^1$ function with
$|\nabla\beta(X)|\in L^{\infty}(\R^N)$, and $R\in
\left(L^2(0,T;L^{2}(\Omega))\right)^N,$ $u\in
L^2(0,T;H^1_0(\Omega))$ satisfy
\begin{equation}\label{R-equation}
\frac{\partial}{\partial_t}R+\mathrm{div}(u \otimes R)=0,\;\; R|_{t=0}=R_0(x)
\end{equation} in the distribution sense.
Then we have
\begin{equation}
\label{renormalized} (\beta(R))_t+\mathrm{div}(\beta(R)
u)+[\nabla\beta(R)\cdot R-\beta(R)]\mathrm{div}u=0
\end{equation}
in the distribution sense. Moreover, if $R\in L^{\infty}(0,T;L^{\gamma}(\Omega))$ for $\gamma>1$, then
$$R\in C([0,T];L^1(\Omega)),$$ and so
$$\int_{\Omega}\beta(R)\,dx(t)=\int_{\Omega}\beta(R_0)\,dx-\int_0^t\int_{\Omega}[\nabla\beta(R)\cdot
R-\beta(R)]\mathrm{div}u\,dx\,dt.$$
\end{lemma}
\begin{remark}
If $N=1$, it is the result of  Feireisl \cite{Feireisl2}.
\end{remark}

 To prove Lemma \ref{main lemma}, we shall rely on the following lemma
which was called the commutator lemma.
\begin{lemma}\label{le:2.3} \cite{Lions}. There exists $C>0$ such that
for any $\rho\in L^2(\mathbb R^d)$ and $u\in H^1(\mathbb R^d)$,
$$
\left\|\eta_\sigma\ast {\rm{div}}(\rho
u)-{\rm{div}}(u(\rho\ast\eta_\sigma))\right\|_{L^1(\mathbb
R^d)}\leq C\|u\|_{H^1(R^d)}\|\rho\|_{L^2(\mathbb R^d)}.
$$
In addition,
$$
\eta_\sigma\ast {\rm{div}}(\rho
u)-{\rm{div}}(u(\rho\ast\eta_\sigma))\rightarrow0 \ in\
L^1(\mathbb R^d), \ as\ \sigma\rightarrow0,$$ where
$\eta_{\sigma}=\frac{1}{\sigma^{d}}\eta(\frac{x}{\sigma})$, and
$\eta(x)\geq 0$ is a smooth even function compactly supported in
the space ball of radius $1$, and with integral equal to $1$.
\end{lemma}
\textbf{Proof of Lemma \ref{main lemma}.}  Here, we devote the  proof of Lemma \ref{main lemma}. The first two steps are similar to the work of \cite{Feireisl}.
\begin{flushleft}
    Step 1: Proof of \eqref{renormalized}.
\end{flushleft}
 Applying the
regularizing operator $f\longmapsto f\ast \eta_{\sigma}$ to both
sides of \eqref{R-equation}, we obtain
\begin{equation}
\label{R equation with sigma}
(R_{\sigma})_t+\mathrm{div}(u\otimes R_{\sigma})=S^{\sigma},
\end{equation}
 almost everywhere on $O\subset\bar{O}\subset(0,T)\times\Omega$ provided $\sigma>0$ small enough, where
 $$S^{\sigma}=\mathrm{div}(u\otimes R_{\sigma})-(\rm{div}(u\otimes R ))_{\sigma},$$
  and $f_{\sigma}(x)=f\ast \eta_{\sigma}.$
  Thanks to Lemma \ref{le:2.3}, we conclude that $$S^{\sigma}\to0 \;\;\text{ in } L^1(O)\;\;\text{ as } \sigma\to 0.$$
 Equation \eqref{R equation with sigma} multiplied by $\nabla\beta(R)$, where $\beta$ is a $C^{1}$ function, gives us
 $$[\beta(R_{\sigma})]_t+\rm{div}[\beta(R_{\sigma})u]+[\nabla\beta(R_{\sigma})\cdot R_{\sigma}-\beta(R_{\sigma})]\rm{div} u=\nabla\beta(R_{\sigma})\cdot S^{\sigma}.$$
This yields \eqref{renormalized} by letting $\sigma\to 0$.

\begin{flushleft}
    Step 2: Continuity of $R$ in the strong topology.
\end{flushleft}
By \eqref{renormalized}, we have
\begin{equation}
\label{renormalized cutoff}
\frac{\partial}{\partial_t}T_K(R)+\mathrm{div}(T_K(R) u) +(\nabla T_K(R)\cdot R-T_K(R))\mathrm{div }u=0\quad\text{ in }\mathfrak{D}'((0,T)\times\Omega),
\end{equation}
 where $T_k(R)$ is a cutoff function verifying $$T_K(R)=\widetilde{T_K}(|R|),\;\;\text{ and } \widetilde{T_K}(z)=K T(\frac{z}{K}),$$
 and $T(z)=z\;\text{ for any } z\in[0,1],\;\;\text{ and it is concave on }\;[0,\infty),\;\;T(z)=2\;\text{ for any} z\geq 3,$
 and $T$ is a $C^{\infty}$ function. We conclude that $T_K(R)$ is bounded in $C([0,T];L^2_{weak}(\Omega))$ due to $R\in L^{\infty}(0,T;L^2(\Omega)).$
Thanks to \eqref{renormalized cutoff},
we have \begin{equation}
\label{weak space}
T_K(R)\;\;\text{ belong to }\, C([0,T];L^{\gamma}_{\text{weak}}(\Omega))
\end{equation}
for any $K\geq 1.$

Applying the same argument
as in step 1 for \eqref{renormalized cutoff}, we have
\begin{equation}
\label{smooth cutoff equation}
\frac{\partial}{\partial_t}[T_K(R)]_{\sigma}+\mathrm{div}([T_K(R)]_{\sigma} u)=A_K^{\sigma}\;\;\; \text{ a.e. on }(0,T)\times U,
\end{equation}
 where $U$ is a compact subset of $\Omega$. Thanks to Lemma \ref{le:2.3}, $A_K^{\sigma}$ is bounded in $L^2(0,T;L^1(\Omega)).$
 Meanwhile, using $2[T_K(R)]_{\sigma}$
to multiply \eqref{smooth cutoff equation},
we have
$$
\frac{\partial}{\partial_t}([T_K(R)]_{\sigma})^2+\mathrm{div}(([T_K(R)]_{\sigma})^2 u)+([T_K(R)]_{\sigma})^2 \mathrm{div} u=2[T_K(R)]_{\sigma} A_K^{\sigma}.$$
Thus, for any test function $\eta\in \mathcal{D}(\Omega),$ the family of functions with respect to $\sigma$ for fixed $K$
$$t\longmapsto \int_{\Omega}([T_K(R)]_{\sigma})^2(t,x)\eta(x)\,dx,\;\sigma>0\;\text{ is precompact in } C[0,T].$$
Note that $[T_K(R)]_{\sigma}\to [T_K(R)]$ in $L^2(\Omega)$ for any $t\in [0,T]$ as $\sigma\to 0$, we obtain
$$t\longmapsto \int_{\Omega}([T_K(R)])^2(t,x)\eta(x)\,dx \text{ is in } C[0,T]$$
for any fixed $\eta(x).$ Thus, $T_K(R)\in C([0,T];L^2(\Omega))$ for any fixed $K\geq 1.$
It allows us to  have $$R\in (C([0,T];L^1(\Omega)))^N,$$
 thanks to \eqref{weak space}.
\begin{flushleft}
    Step 3: Final inequality.
\end{flushleft}
Taking integration on \eqref{renormalized} with respect
to $t$, we have
  $$\int_{\Omega}\beta(R(t))\,dx=\int_{\Omega}\beta(R(s))\,dx-\int_s^t\int_{\Omega}[\nabla\beta(R)\cdot R-\beta(R)]\mathrm{div}u\,dx\,dt,$$
  where $0<s<t<T.$ Thanks to $$R\in (C([0,T];L^1(\Omega)))^N.$$

  Letting $s\to 0$, thus we have

  $$\int_{\Omega}\beta(R(t))\,dx=\int_{\Omega}\beta(R_0)\,dx-\int_0^t\int_{\Omega}[\nabla\beta(R)\cdot R-\beta(R)]\mathrm{div}u\,dx\,dt$$
  for any $0\leq t\leq T.$
\endpf

\bigskip
With above lemmas in hand, we are ready to show Theorem \ref{main 2-le:important}.\\
\newline
\textbf{Proof of Theorem \ref{main 2-le:important}. } Up to a
subsequence, \begin{equation}\label{weak convergence for general}
\rho_K\to\rho,\;n_K\to n\;\text{ weakly in }
L^{\infty}(0,T;L^2(\Omega)),\quad
 u_K\to u \text{ weakly in } L^2(0,T;H^1_0(\Omega)),
\end{equation} as $K\rightarrow\infty$.
Passing to the limit as $K\to\infty$ in \eqref{the parabolic for
density} and \eqref{second parabolic equation for n} respectively,
we have
\begin{equation*}
\label{parabolic for density}
\rho_t+\mathrm{div}(\rho u)=0,\quad\;\;\rho|_{t=0}=\rho_0,
\end{equation*}
and
\begin{equation*}
\label{parabolic equation for n} n_t+\mathrm{div}(n
u)=0,\quad\;\;n|_{t=0}=n_0.
\end{equation*}
Using Lemma \ref{main lemma} with $R=(b,d)$ and $u=u$,
$\beta_{\sigma}(b,d)=\frac{b^2}{d+\sigma}$, note that
$$\nabla\beta_\sigma(R)\cdot
R-\beta_\sigma(R)=\sigma\frac{b^2}{(d+\sigma)^2},$$ one obtains
$$\int_{\Omega}\frac{b(t,x)^2}{d(t,x)+\sigma}\,dx=\int_{\Omega}\frac{b_0^2}{d_0+\sigma}\,dx-
\sigma\int_0^t\int_{\Omega}\frac{b^2}{(d+\sigma)^2}\mathrm{div}u\,dx\,dt,$$
for almost everywhere $t\in[0,T].$\\
Note that $$\frac{b^2}{d+\sigma}\leq \frac{b^2}{d}\;\;\text{ and
}\frac{b_0^2}{d_0+\sigma}\leq \frac{b_0^2}{d_0},$$ by the
dominated convergence theorem, we obtain the following equality by
letting $\sigma$ goes to zero,
\begin{equation}
\label{key estimate-equality}
\int_{\Omega}\frac{b(t,x)^2}{d(t,x)}\,dx=\int_{\Omega}\frac{b_0^2}{d_0}\,dx
\end{equation}
for almost everywhere $t\in[0,T].$

By
\eqref{AAAinitial condition for n2 over density} and \eqref{key estimate-equality}, we find
\begin{equation}
\label{key condition} \int_\Omega \frac{b_K(t,x)^2}
{d_K(t,x)}\,dx=\int_{\Omega}b_Ka_K\,dx\le \int_\Omega
\frac{b(t,x)^2} {d(t,x)}\,dx=\int_{\Omega}b a\,dx
\end{equation}
for almost everywhere $t\in[0,T].$

Thanks to \eqref{weak convergence for general},  setting
$(b_K,d_K)=(h_K,g_K)$ for Lemma \ref{2-le:important}, one obtains
\eqref{key bound} for any $s>1.$
\endpf
\newline

\section{Faedo-Galerkin approach}
\setcounter{equation}{0} \setcounter{theorem}{0}
In this section, we construct a global weak solution $(\rho,n,u)$ to the following approximation  \eqref{a-equation}-\eqref{a-boundary} with a finite energy.
Motivated by the work of \cite{Feireisl}, we propose the following approximation system
\begin{equation}\label{a-equation}
    \left\{
    \begin{array}{l}
        n_t+\mathrm{div}(n u)=\epsilon\Delta n,\\[2mm]
        \rho_t+\mathrm{div}(\rho u)=\epsilon\Delta \rho,\\[2mm]
 \big[(\rho+n) u\big]_t+\mathrm{div}\big[(\rho+n) u\otimes u\big]
 +\nabla(n^\alpha+\rho^\gamma)+\delta\nabla(\rho+n)^\beta+\epsilon\nabla u\cdot\nabla(\rho+n)\\ =
            \mu \Delta u + (\mu+\lambda)\nabla \mathrm{div}u
    \end{array}
    \right.
    \end{equation}
    on  $\Omega\times (0, \infty)$, with initial and boundary condition
\be\label{a-initial} \big(\rho, n,(\rho+n)
u\big)|_{t=0}=(\rho_{0,\delta},n_{0,\delta},M_{0,\delta})\
\mathrm{on}\ \overline{\Omega}, \ee \be\label{a-boundary}
(\frac{\partial\rho}{\partial \nu},\frac{\partial n}{\partial\nu},
u)|_{\partial\Omega}=0, \ee where $\epsilon,\delta>0$,
$\beta>\max\{\alpha,\gamma\}$,
$M_{0,\delta}=(\rho_{0,\delta}+n_{0,\delta})u_{0,\delta}$ and
$n_{0,\delta},\rho_{0,\delta}\in C^3(\overline{\Omega})$,
$u_{0,\delta}\in C_0^3(\Omega)$ satisfying
\be\label{appinitial}\begin{cases} 0<\delta\le
\rho_{0,\delta},n_{0,\delta}\le \delta^{-\frac{1}{2\beta}},\quad
 (\frac{\partial n_{0,\delta}}{\partial
\nu}, \frac{\partial
\rho_{0,\delta}}{\partial \nu})|_{\partial\Omega}=0,\\[2mm]
\lim\limits_{\delta\rightarrow0}\left(\|\rho_{0,\delta}-\rho_0\|_{L^\gamma(\Omega)}+\|n_{0,\delta}-n_0\|_{L^\alpha(\Omega)}\right)=0,\\[2mm]
u_{0\delta}=\frac{\varphi_\delta}{\sqrt{\rho_{0,\delta}+n_{0,\delta}}}\eta_\delta*(\frac{M_0}{\sqrt{\rho_0+n_0}}),
\\[2mm]
\sqrt{\rho_{0,\delta}+n_{0,\delta}}u_{0,\delta}\rightarrow\frac{M_0}{\sqrt{\rho_0+n_0}}\quad
\mathrm{in}\
L^2(\Omega) \ \mathrm{as}\ \delta\rightarrow0,\\[2mm]
m_{0,\delta}\rightarrow M_0\quad \mathrm{in}\ L^1(\Omega) \
\mathrm{as}\ \delta\rightarrow0,\\[2mm]
 \frac{1}{c_0}\rho_{0,\delta}\le n_{0,\delta}\le
c_0\rho_{0,\delta}\quad \mathrm{if}\,\, \frac{1}{c_0}\rho_{0}\le
n_{0}\le c_0\rho_{0},
\end{cases} \ee where $\delta\in(0,1)$, $\eta$ is the standard mollifier, $\varphi_\delta\in C_0^\infty(\Omega)$, $0\le\varphi_\delta\le1$ on
$\overline{\Omega}$ and $\varphi_\delta\equiv 1$ on
$\big\{x\in\Omega| \mathrm{dist}(x,\partial\Omega)>\delta\big\}$.


We are able to use Faedo-Galerkin approach to construct a
global weak solution to (\ref{a-equation}), (\ref{a-initial}) and
(\ref{a-boundary}). To begin with, we consider a sequence of
finite dimensional spaces \bex
X_k=[span\{\psi_j\}_{j=1}^k]^3,\quad
k\in\{1,2,3,\cdot\cdot\cdot\}, \eex where
$\{\psi_i\}_{i=1}^\infty$ is the set of the eigenfunctions of the
Laplacian: \bex\begin{cases} -\Delta\psi_i=\lambda_i\psi_i \quad \mathrm{on}\ \Omega,\\[2mm]
\psi_i|_{\partial\Omega}=0.
\end{cases}\eex

For any given $\epsilon,\delta>0$, we shall look for the
approximate solution $u_k\in C([0,T];X_k)$ (for any fixed $T>0$)
given by the following form: \be\label{1-apu}\begin{split}
&\int_\Omega(\rho_k+n_k) u_k(t)\cdot\psi\,dx-\int_\Omega
m_{0,\delta}\cdot\psi\,dx=\int_0^t\int_\Omega\big[\mu \Delta u_k +
(\mu+\lambda)\nabla
\mathrm{div}u_k\big]\cdot\psi\,dx\,ds\\&-\int_0^t\int_\Omega\Big[\mathrm{div}\big[(\rho_k+n_k)
u_k\otimes u_k\big]
 +\nabla(n_k^\alpha+\rho_k^\gamma)+\delta\nabla(\rho_k+n_k)^\beta+\epsilon\nabla
 u_k\cdot\nabla(\rho_k+n_k)\Big]\cdot\psi\,dx\,ds
\end{split}
\ee for $t\in[0, T]$ and $\psi\in X_k$, where $\rho_k=\rho_k(u_k)$
and $n_k=n_k(u_k)$ satisfying \be\label{1-apnrho}\begin{cases}
\partial_tn_k+\mathrm{div}(n_k u_k)=\epsilon\Delta n_k,\\[2mm]
        \partial_t\rho_k+\mathrm{div}(\rho_k u_k)=\epsilon\Delta \rho_k,\\[2mm]
        n_k|_{t=0}=n_{0,\delta},\quad
        \rho_k|_{t=0}=\rho_{0,\delta},\\[2mm]
(\frac{\partial\rho_k}{\partial \nu},\frac{\partial
n_k}{\partial\nu})|_{\partial\Omega}=0.
\end{cases}
\ee


Due to Lemmas 2.1 and 2.2 in \cite{Feireisl}, the problem
(\ref{1-apu}) can be solved on a short time interval $[0, T_k]$
for $T_k\le T$ by a standard fixed point theorem on the Banach
space $C([0, T_k]; X_k)$.
To show that $T_k=T$, we need the uniform estimates resulting  from the following energy equality
 \be\label{0-energy}
\begin{split} &\frac{d}{dt}\int_\Omega\Big[\frac{(\rho_k+n_k)|u_k|^2}{2}+G_\alpha(n_k)+\frac{1}{\gamma-1}\rho_k^\gamma+\frac{\delta}{\beta-1}(\rho_k+n_k)^\beta\Big]\,dx
\\&+\int_\Omega\Big[\mu|\nabla u_k|^2
+(\mu+\lambda)|\mathrm{div}u_k|^2\Big]\,dx\\&+\int_\Omega\Big[\epsilon\alpha
n_k^{\alpha-2}|\nabla n_k|^2+\epsilon\gamma
 \rho_k^{\gamma-2}|\nabla\rho_k|^2 +\epsilon\beta\delta(\rho_k+n_k)^{\beta-2}|\nabla(\rho_k+n_k)|^2\Big]\,dx
=0,\ \mathrm{on}\ (0,T_k),
\end{split}
\ee  where \bex
G_\alpha(n_k)=\left\{\begin{array}{l} n_k\ln
n_k-n_k+1, \ \ \ {\rm for} \ \ \ \alpha=1, \\
[3mm] \frac{n_k^\alpha}{\alpha-1}, \ \ \ {\rm for}\ \ \ \alpha>1.
\end{array}
\right.\eex
This could be done by  differentiating (\ref{1-apu}) with respect to time, taking
$\psi=u_k(t)$ and using (\ref{1-apnrho}). We refer the readers to \cite{Feireisl} for more
details. Thus, we obtain a solution $(\rho_k, n_k, u_k)$ to \eqref{1-apu}-\eqref{1-apnrho} globally in time $t$.

The next step
is to  pass the limit of $(\rho_k, n_k, u_k)$ as
$k\rightarrow\infty$.
Following the same arguments of Section 2.3 of \cite{Feireisl},  energy equality \eqref{0-energy} gives us the following bounds

\be\label{1-0} 0<\frac{1}{c_k}\le \rho_k(x,t),n_k(x,t)\le c_k\
\mathrm{for}\ \mathrm{a.e.} (x,t)\in\Omega\times (0,T), \ee

\be\label{1-1}
\sup\limits_{t\in[0,T]}\|\rho_k(t)\|_{L^\gamma(\Omega)}^\gamma\le
C(\rho_0,n_0,M_0),\ee

\be\label{1-1+1}
\sup\limits_{t\in[0,T]}\|n_k(t)\|_{L^\alpha(\Omega)}^\alpha\le
C(\rho_0,n_0,M_0)\quad \mathrm{for}\ \alpha\ge1,\ee

\be\label{1-2}\delta\sup\limits_{t\in[0,T]}\|\rho_k(t)+n_k(t)\|_{L^\beta(\Omega)}^\beta\le
C(\rho_0,n_0,M_0),\ee

\be\label{1-3}\sup\limits_{t\in[0,T]}\|\sqrt{\rho_k+n_k}(t)u_k(t)\|_{L^2(\Omega)}^2\le
C(\rho_0,n_0,M_0),\ee

 \be\label{1-4}
\int_0^T\|u_k(t)\|_{H^1_0(\Omega)}^2\,dt\le C(\rho_0,n_0,M_0),\ee

\be\label{1-5}
\epsilon\int_0^T\big(\|\nabla\rho_k(t)\|_{L^2(\Omega)}^2+\|\nabla
n_k(t)\|_{L^2(\Omega)}^2\big)\,dt\le
C(\beta,\delta,\rho_0,n_0,M_0),\ee and
 \be\label{1-6}
\|\rho_k+n_k\|_{L^{\beta+1}(Q_T)}\le
C(\epsilon,\beta,\delta,\rho_0,n_0,M_0), \ee where
$Q_T=\Omega\times (0,T)$ and $\beta\ge4$.
\\
\newline
We are able to control $n_k$ by $\rho_k$ if some additional initial data is given in \eqref{additional initial data}.
\begin{lemma}
\label{lemma bound on n} If $(\rho_k,n_k,u_k)$ is a solution to
(\ref{1-apu}) and (\ref{1-apnrho}) with the initial data
satisfying \begin{equation}
\label{additional initial data}
\frac{1}{c_0}\rho_0\le n_0\le c_0\rho_0
\end{equation} on $\Omega$,then the following inequality holds \be\label{1-6+1}
\frac{1}{c_0}\rho_k(x,t)\le n_k(x,t)\le c_0\rho_k(x,t)\ee for a.e.
$(x,t)\in Q_T$.
\end{lemma}

\pf It is easy to check that $n_k-c_0\rho_k$ is a solution of the
following parabolic equation
\begin{equation*}\begin{cases}
(n_k-c_0\rho_k)_t+\mathrm{div}\big[(n_k-c_0\rho_k)u_k\big]=\epsilon\Delta(n_k-c_0\rho_k),\\[2mm]
(n_k-c_0\rho_k)|_{t=0}=n_{0,\delta}-c_0\rho_{0,\delta},\\[2mm]
\nabla(n_k-c_0\rho_k)\cdot\nu|_{\partial\Omega}=0.
\end{cases}
\end{equation*}
The right inequality of (\ref{1-6+1}) can be obtained by applying
the maximum principle on it. Similarly, we obtain the left
inequality of (\ref{1-6+1}).

\endpf

If the initial data satisfies \eqref{additional initial data}  and with (\ref{1-1}), (\ref{1-1+1}), and
(\ref{1-6+1}), we have \be\label{1-7}
\sup\limits_{t\in[0,T]}\Big(\|\rho_k(t)\|_{L^{\alpha_1}(\Omega)}^{\alpha_1}+\|n_k(t)\|_{L^{\alpha_1}(\Omega)}^{\alpha_1}\Big)\le
C(\rho_0,n_0,M_0),\ee where $\alpha_1=\max\{\alpha,\gamma\}$.

Relying on the above uniform estimates, i.e.,
(\ref{1-1})-(\ref{1-6+1}) and (\ref{1-7}), and the Aubin-Lions
lemma, we are able to recover the global solution to the
approximation system (\ref{a-equation})-(\ref{a-boundary}) by
passing to the limit for $(\rho_k, n_k, u_k)$ as
$k\rightarrow\infty$. We have the following Proposition on the
weak solutions of the approximation  (\ref{a-equation}),
(\ref{a-initial}) and (\ref{a-boundary}).
\begin{proposition}\label{0-le:aweak solution}
Suppose $\beta>\max\{4,\alpha,\gamma\}$. For any given
$\epsilon,\delta>0$, there exists a global weak solution
$(\rho,n,u)$ to (\ref{a-equation}), (\ref{a-initial}) and
(\ref{a-boundary}) such that for any given $T>0$, the following
estimates
 \be\label{1-r1}
\sup\limits_{t\in[0,T]}\|\rho(t)\|_{L^\gamma(\Omega)}^\gamma\le
C(\rho_0,n_0,M_0),\ee \be\label{1-r1+1}
\sup\limits_{t\in[0,T]}\|n(t)\|_{L^\alpha(\Omega)}^\alpha\le
C(\rho_0,n_0,M_0),\ee
\be\label{1-r2}\delta\sup\limits_{t\in[0,T]}\|(\rho(t),n(t))\|_{L^\beta(\Omega)}^\beta\le
C(\rho_0,n_0,M_0),\ee
\be\label{1-r3}\sup\limits_{t\in[0,T]}\|\sqrt{\rho+n}(t)u(t)\|_{L^2(\Omega)}^2\le
C(\rho_0,n_0,M_0),\ee \be\label{1-r4}
\int_0^T\|u(t)\|_{H^1_0(\Omega)}^2\,dt\le C(\rho_0,n_0,M_0),\ee
\be\label{1-r5} \epsilon\int_0^T\|(\nabla\rho(t),\nabla
n(t))\|_{L^2(\Omega)}^2\,dt\le C(\beta,\delta,\rho_0,n_0,M_0),\ee
and \be\label{1-r6} \|(\rho(t),n(t))\|_{L^{\beta+1}(Q_T)}\le
C(\epsilon,\beta,\delta,\rho_0,n_0,M_0) \ee hold, where the norm
$\|(\cdot,\cdot)\|$ denotes $\|\cdot\| + \|\cdot\|$, and
$\rho,n\ge0$ a.e. on $Q_T$.

In addition, if the initial data satisfy $\frac{1}{c_0}\rho_0\le
n_0\le c_0\rho_0$ on $\Omega$,
then\be\label{nrho}\begin{cases}\frac{1}{c_0}\rho\le n\le
c_0\rho\quad \mathrm{a.e.}\,\, \mathrm{on}\,\, \Omega\times
(0,T),\\[2mm]
\sup\limits_{t\in[0,T]}\|(\rho,n)(t)\|_{L^{\alpha_1}(\Omega)}^{\alpha_1}\le
C(\rho_0,n_0,M_0),\end{cases}\ee where
$\alpha_1=\max\{\alpha,\gamma\}$. Finally, there exists $r>1$ such
that $\rho_t, n_t, \nabla^2\rho, \nabla^2 n\in L^r(Q_T)$ and the
equations (\ref{a-equation})$_1$ and (\ref{a-equation})$_2$ are
satisfied a.e. on $Q_T$.
\end{proposition}

\begin{remark}
The solution $(\rho,n,u)$ stated in Proposition \ref{0-le:aweak
solution} actually depends on $\epsilon$ and $\delta$. We omit the
dependence in the solution itself for brevity.
\end{remark}

\section{The vanishing viscosity limit $\epsilon\rightarrow0^+$}
\setcounter{equation}{0} \setcounter{theorem}{0}

The goal of this section is to pass to the limit of
$(\rho_{\epsilon},n_{\epsilon},u_{\epsilon})$ as $\epsilon$ goes
to zero. To vanish $\epsilon$, the uniform estimates are needed.
Compared to the work of \cite{Feireisl}, the pressure law involves
 two variables, which bring new difficulty-possible oscillation
of $\rho^{\gamma}+n^{\alpha}$. The uniform estimates resulting
from the energy inequality in Proposition \ref{0-le:aweak
solution} and Lemma \ref{le:h-inofrho} are not enough to handle
the weak limit of such a pressure. In Section \ref{s4.1}, we pass
to the limits for the weak solution constructed in Proposition
\ref{0-le:aweak solution} as $\epsilon$ goes to zero by standard
compactness argument. In Section \ref{s4.2}, we will focus on the
weak limit of the pressure and the strong convergence of
$\rho_{\epsilon}$ and $n_{\epsilon}$. In this section, let C
denote a generic positive constant depending on the initial data
and $\delta$ but independent of $\epsilon$.

\subsection{Passing to the limit as
$\epsilon\rightarrow0^+$}\label{s4.1}

The uniform estimates resulting from (\ref{1-r1}), (\ref{1-r1+1}),
(\ref{1-r2}), and (\ref{nrho}) are not enough to obtain the
convergence of the pressure term
$\rho_{\epsilon}^{\gamma}+n_{\epsilon}^{\alpha}$. Thus we need to
obtain higher integrability estimates of the pressure term
uniformly for $\epsilon$.

First, following the same argument in \cite{Feireisl}, we are able to get the following estimate in Lemma \ref{le:h-inofrho}.
\begin{lemma}\label{le:h-inofrho}
Let $(\rho,n,u)$ be the solution stated in Proposition \ref{0-le:aweak
solution}, then \bex \int_0^T\int_\Omega
(n^{\alpha+1}+\rho^{\gamma+1}+\delta\rho^{\beta+1}+\delta
n^{\beta+1})\,dx\,dt\le C\eex for $\beta>4$.
\end{lemma}

\bigskip

 In this step, we fix $\delta>0$ and shall let
$\epsilon\rightarrow0^+$. Then the solution $(\rho,n,u)$
constructed in Proposition \ref{0-le:aweak solution} is naturally
dressed in the subscript ``$\epsilon$", i.e.,
$(\rho_\epsilon,n_\epsilon, u_\epsilon)$.

With (\ref{1-r1})-(\ref{1-r5}), and Lemma \ref{le:h-inofrho},
letting $\epsilon\rightarrow0^+$ (take the subsequence if
necessary), we have \be\label{2-lim}
\begin{cases}
(\rho_\epsilon,n_\epsilon)\rightarrow (\rho,n)\ \mathrm{in}\
C([0,T];
L_{weak}^\beta(\Omega))\ \mathrm{and}\ \mathrm{weakly}\ \mathrm{in}\ L^{\beta+1}(Q_T)\ \mathrm{as}\
\epsilon\rightarrow0^+,\\[2mm]
(\epsilon\Delta\rho_\epsilon, \epsilon\Delta
n_\epsilon)\rightarrow 0\ \mathrm{weakly}\ \mathrm{in}\ L^2(0,T;
H^{-1}(\Omega))\ \mathrm{as}\ \epsilon\rightarrow0^+,\\[2mm]
u_\epsilon\rightarrow u\ \mathrm{weakly}\ \mathrm{in}\ L^2(0,T;
H_0^1(\Omega))\ \mathrm{as}\ \epsilon\rightarrow0^+,\\[2mm]
(\rho_\epsilon+n_\epsilon) u_\epsilon\rightarrow (\rho+n) u \
\mathrm{in}\ C([0,T]; L_{weak}^\frac{2\gamma}{\gamma+1})\cap
C([0,T];H^{-1}(\Omega))\ \mathrm{as}\ \epsilon\rightarrow0^+,\\[2mm]
(\rho_\epsilon u_\epsilon,n_\epsilon u_\epsilon)\rightarrow (\rho
u, nu) \ \mathrm{in}\ \mathcal{D}^\prime(Q_T)\ \mathrm{as}\
\epsilon\rightarrow0^+,\\[2mm]
(\rho_\epsilon+n_\epsilon) u_\epsilon\otimes u_\epsilon\rightarrow
(\rho+n) u\otimes u\ \mathrm{in}\ \mathcal{D}^\prime(Q_T)\
\mathrm{as}\
\epsilon\rightarrow0^+,\\[2mm]
n_\epsilon^\alpha+\rho_\epsilon^\gamma+\delta(\rho_\epsilon+n_\epsilon)^\beta\rightarrow
\overline{n^\alpha+\rho^\gamma+\delta(\rho+n)^\beta}\
\mathrm{weakly}\ \mathrm{in}\ L^\frac{\beta+1}{\beta}(Q_T)\
\mathrm{as}\
\epsilon\rightarrow0^+,\\[2mm]
\epsilon\nabla
u_\epsilon\cdot\nabla(\rho_\epsilon+n_\epsilon)\rightarrow 0\
\mathrm{in}\ L^1(Q_T)\ \mathrm{as}\ \epsilon\rightarrow0^+,
\end{cases}
\ee and $\rho,n\ge0$.

 By virtue of
(\ref{nrho}) and (\ref{2-lim})$_1$, if $\frac{1}{c_0}\rho_0\le
n_0\le c_0\rho_0$, we have \bex\frac{1}{c_0}\rho_\epsilon(x,t)\le
n_{\epsilon}(x,t)\le c_0\rho_{\epsilon}(x,t)\ \mathrm{and}\
\frac{1}{c_0}\rho(x,t)\le n(x,t)\le c_0 \rho(x,t)\quad
\mathrm{for}\ \mathrm{a.e.}\ (x,t)\in Q_T.\eex

With (\ref{2-lim})$_1$ and (\ref{2-lim})$_4$, we get $$\big(\rho,
n,(\rho+n)
u\big)|_{t=0}=(\rho_{0,\delta},n_{0,\delta},M_{0,\delta}).$$

In summary, the limit $(\rho, n, u)$ solves the following system
in the sense of distribution on $Q_T$ for any $T>0$:
 \begin{equation}\label{a2-equation}
    \left\{
    \begin{array}{l}
        n_t+\mathrm{div}(n u)=0,\\[2mm]
        \rho_t+\mathrm{div}(\rho u)=0,\\[2mm]
 \big[(\rho+n) u\big]_t+\mathrm{div}\big[(\rho+n) u\otimes u\big]
 +\nabla(\overline{n^\alpha+\rho^\gamma+\delta(\rho+n)^\beta}) =
            \mu \Delta u + (\mu+\lambda)\nabla \mathrm{div}u
    \end{array}
    \right.
    \end{equation}
    with initial and boundary condition
\be\label{a2-initial} \big(\rho, n,(\rho+n)
u\big)|_{t=0}=(\rho_{0,\delta},n_{0,\delta},M_{0,\delta}), \ee
\be\label{a2-boundary} u|_{\partial\Omega}=0, \ee
where $\overline{f(t,x)}$ denotes the weak limit of $f_{\epsilon}(t,x)$ as $\epsilon\to0.$

To this end, we have to show that $\overline{n^\alpha+\rho^\gamma+\delta(\rho+n)^\beta}=n^\alpha+\rho^\gamma+\delta(\rho+n)^\beta$, which is a nonlinear two-variable function in term of $\rho$ and $n$. It seems that the argument in \cite{Feireisl} fails here due to the difficulty resulting from the new variable $n$. New ideas are necessary to handle this weak limit. We are going to focus on this issue next subsection.

\subsection{The weak limit of pressure}\label{s4.2}
The main task of this subsection is to
handle the possible oscillation for the pressure
$n_{\epsilon}^\alpha+\rho_{\epsilon}^\gamma+\delta(\rho_{\epsilon}+n_{\epsilon})^\beta$.
To achieve this goal, we have to show the strong convergence of
$\rho_{\epsilon}$ and $n_{\epsilon}$. It allows us to have the
following Proposition on the weak limit of pressure.
\begin{proposition}
\label{key pro}
$$
\overline{n^\alpha+\rho^\gamma+\delta(\rho+n)^\beta}=n^\alpha+\rho^\gamma+\delta(\rho+n)^\beta
$$a.e. on $Q_T$.
\end{proposition}

To prove this proposition, we shall rely on the following lemmas. The first one is on the effective viscous flux associated with pressure involving two variables. In particular, let
 \bex\begin{split}
H_\epsilon:=&n_\epsilon^\alpha+\rho_\epsilon^\gamma+\delta(\rho_\epsilon+n_\epsilon)^\beta-(2\mu+\lambda)\mathrm{div}u_\epsilon,
\\
H:=&\overline{n^\alpha+\rho^\gamma+\delta(\rho+n)^\beta}-(2\mu+\lambda)\mathrm{div}u,
\end{split}
\eex
then we will have the following lemma. The proof is very similar to the work of \cite{Feireisl}.
\begin{lemma}\label{2-le:3.7} Let
$(\rho_\epsilon,n_\epsilon,u_\epsilon)$ be the solution stated in
Lemma \ref{0-le:aweak solution}, and $(\rho, n,u)$ be the limit in
the sense of (\ref{2-lim}), then \be\label{2-8}
\lim\limits_{\epsilon\rightarrow0^+}\int_0^T\psi\int_\Omega\phi
H_\epsilon(\rho_\epsilon+n_\epsilon)\,dx\,dt=\int_0^T\psi\int_\Omega\phi
H(\rho+n)\,dx\,dt, \ee for any $\psi\in  C_0^\infty(0,T)$ and
$\phi\in
 C_0^\infty(\Omega)$.
\end{lemma}

The key idea of proving Proposition \ref{key pro} is to rewrite the terms related  pressure  as follows
\bex
\begin{split}
n_\epsilon^\alpha+\rho_\epsilon^\gamma=&A_\epsilon^\alpha
d_\epsilon^\alpha+B_\epsilon^\gamma d_\epsilon^\gamma =A^\alpha
d_\epsilon^\alpha+B^\gamma
d_\epsilon^\gamma+(A_\epsilon^\alpha-A^\alpha)
d_\epsilon^\alpha+(B_\epsilon^\gamma-B^\gamma) d_\epsilon^\gamma,\\
n_\epsilon+\rho_\epsilon=&(A_\epsilon+B_\epsilon) d_\epsilon
=(A+B)d_\epsilon +(A_\epsilon-A+B_\epsilon-B) d_\epsilon,
\end{split}
\eex where $d_\epsilon=\rho_\epsilon+n_\epsilon$, $d=\rho+n$,
$(A_\epsilon,B_\epsilon)=(\frac{n_\epsilon}{d_\epsilon},\frac{\rho_\epsilon}{d_\epsilon})$
if $d_\epsilon\neq0$, $(A,B)=(\frac{n}{d},\frac{\rho}{d})$ if
$d\neq0$, $0\le A_\epsilon,B_\epsilon, A, B\le 1$, and
$(A_\epsilon d_\epsilon,B_\epsilon
d_\epsilon)=(n_\epsilon,\rho_\epsilon), (Ad,Bd)=(n,\rho)$,
$(\rho,n)$ is the limit of $(\rho_{\epsilon},n_{\epsilon})$ in a
suitable weak topology. We are able to apply the ideas in
\cite{Feireisl} to handle the product $A^\alpha
d_\epsilon^\alpha+B^\gamma d_\epsilon^\gamma$ and
$(A+B)d_\epsilon$, because $A$ and $B$ are bounded in
$L^{\infty}(0,T;L^{\infty}(\Omega))$ and they are viewed as the coefficients.
The difficult part is to show that the remainder tends to zero as $\varepsilon$ goes to zero.  Theorem \ref{main 2-le:important} allows us to show
the terms $\Big[(A_\epsilon^\alpha-A^\alpha)
d_\epsilon^\alpha+(B_\epsilon^\gamma-B^\gamma)
d_\epsilon^\gamma\Big](n_\epsilon+\rho_\epsilon)$ and
$\Big(A^\alpha d_\epsilon^\alpha+B^\gamma
d_\epsilon^\gamma\Big)(A_\epsilon-A+B_\epsilon-B) d_\epsilon$
approach to zero as $\epsilon$ goes to zero.

We divide the proof of Proposition \ref{key pro} into several
steps as follows:

\begin{flushleft}
\textbf{    Step 1: Control $\rho_{\epsilon}$ and $n_{\epsilon}$ in $L\log L$.}
\end{flushleft}

The current step of our proof is to control $\rho_{\epsilon}$ and $n_{\epsilon}$ in the space of $L\log L$. This helps us to obtain the strong convergence of $\rho_{\epsilon}$ and $n_{\epsilon}$. We give our control in the following lemma.

\begin{lemma}
Let $(\rho_\epsilon,n_\epsilon)$ be the solution stated in Proposition
\ref{0-le:aweak solution}, and $(\rho, n)$ be the limit in the
sense of (\ref{2-lim}), then\be\label{2-17}
\begin{split}
&\int_\Omega \big[\rho_\epsilon\log \rho_\epsilon-\rho\log
\rho+n_\epsilon\log n_\epsilon-n\log n\big](t)\,dx
 \\ \le&\int_0^t\int_\Omega
(\rho+n) \mathrm{div}u\,dx\,ds-\int_0^t\int_\Omega
(\rho_\epsilon+n_\epsilon) \mathrm{div}u_\epsilon\,dx\,ds
\end{split}
\ee for a.e. $t\in(0,T)$.
\end{lemma}
\pf Since  $n_\epsilon$ and $\rho_\epsilon$ solve
(\ref{a-equation})$_1$ and (\ref{a-equation})$_2$ a.e. on $Q_T$,
respectively, we have \be\label{2-14}
\begin{split}
[\beta_j(f_\epsilon)]_t+\mathrm{div}\big(\beta_j(f_\epsilon)u_\epsilon\big)+\big[\beta_j^\prime(f_\epsilon)f_\epsilon-\beta_j(f_\epsilon)\big]
\mathrm{div}u_\epsilon=\epsilon\Delta \beta_j(f_\epsilon)-\epsilon
\beta_j^{\prime\prime}(f_\epsilon)|\nabla f_\epsilon|^2\
\mathrm{on}\ \mathrm{Q_T},
\end{split}
\ee where $f_\epsilon=\rho_\epsilon,n_\epsilon$, and $\beta_j\in
C^2[0,\infty)$.

Taking $\beta_j(z)=(z+\frac{1}{j})\log(z+\frac{1}{j})$ in
(\ref{2-14}), and integrating it over $\Omega\times(0,t)$ for
$t\in[0,T]$, we have \be\label{bbb}
\begin{split}
&\int_\Omega(f_\epsilon+\frac{1}{j})\log(f_\epsilon+\frac{1}{j})(t)\,dx
+\int_0^t\int_\Omega\big[f_\epsilon-\frac{1}{j}\log(f_\epsilon+\frac{1}{j})\big]
\mathrm{div}u_\epsilon\,dx\,ds\\ \le&
\int_\Omega(f_{0,\epsilon}+\frac{1}{j})\log(f_{0,\epsilon}+\frac{1}{j})\,dx,
\end{split}
\ee where we have used the convexity of $\beta_j$ and the boundary
condition (\ref{a-boundary}). Letting $j\rightarrow\infty$ in
\eqref{bbb}, one obtains  \be\label{2-15}
\begin{split}
\int_\Omega (f_\epsilon\log f_\epsilon)(t)\,dx
+\int_0^t\int_\Omega f_\epsilon \mathrm{div}u_\epsilon\,dx\,ds \le
\int_\Omega f_{0,\delta}\log f_{0,\delta}\,dx,
\end{split}
\ee where $f_\epsilon=\rho_\epsilon,n_\epsilon$ and
$f_{0,\delta}=\rho_{0,\delta},n_{0,\delta}$.

Since the limit $(n,u)$ and $(\rho,u)$ solve
(\ref{a2-equation})$_1$ and (\ref{a2-equation})$_2$ in the sense
of renormalized solutions, we can take $\beta(z)=z\log z$ in
accordance with Remark 1.1 in \cite{Feireisl} or by approximating
the function $z\log z$ by a sequence of such the $\beta(z)$ stated in
Lemma \ref{main lemma}  and then passing to the limit. This allows us to have
\be\label{2-16}
\begin{split}
\int_\Omega (f\log f)(t)\,dx +\int_0^t\int_\Omega f
\mathrm{div}u\,dx\,ds =\int_\Omega f_{0,\delta}\log
f_{0,\delta}\,dx,
\end{split}
\ee where $f=\rho,n$ and
$f_{0,\delta}=\rho_{0,\delta},n_{0,\delta}$.
Thanks to (\ref{2-15}) and (\ref{2-16}), (\ref{2-17}) follows.
\begin{flushleft}
\textbf{    Step 2: Control the right hand side of (\ref{2-17})}
\end{flushleft}
It is crucial to control the right hand side of (\ref{2-17}).
Thanks to Theorem \ref{main 2-le:important}, we  show the following lemma which can help us to finish this step.

\begin{lemma}Let $(\rho_\epsilon,n_\epsilon)$ be the solution stated in Proposition
\ref{0-le:aweak solution}, and $(\rho, n)$ be the limit in the
sense of (\ref{2-lim}), then \bex\int_0^t\psi\int_\Omega\phi
(\rho+n) \overline{n^\alpha+\rho^\gamma}\,dx\,ds\le
\int_0^t\psi\int_\Omega\phi \overline{(\rho+n)
(n^\alpha+\rho^\gamma)}\,dx\,ds\eex for any $t\in[0,T]$ and any
$\psi\in C[0,t]$, $\phi\in C(\overline{\Omega})$ where
$\psi,\phi\ge0$.
\end{lemma}

\pf Note that \bex
\begin{split}
n_\epsilon^\alpha+\rho_\epsilon^\gamma=&A_\epsilon^\alpha
d_\epsilon^\alpha+B_\epsilon^\gamma d_\epsilon^\gamma =A^\alpha
d_\epsilon^\alpha+B^\gamma
d_\epsilon^\gamma+(A_\epsilon^\alpha-A^\alpha)
d_\epsilon^\alpha+(B_\epsilon^\gamma-B^\gamma) d_\epsilon^\gamma,\\
n_\epsilon+\rho_\epsilon=&(A_\epsilon+B_\epsilon) d_\epsilon
=(A+B)d_\epsilon +(A_\epsilon-A+B_\epsilon-B) d_\epsilon,
\end{split}
\eex where  $d_\epsilon=\rho_\epsilon+n_\epsilon$, $d=\rho+n$,
$(A_\epsilon,B_\epsilon)=(\frac{n_\epsilon}{d_\epsilon},\frac{\rho_\epsilon}{d_\epsilon})$
if $d_\epsilon\neq0$, $(A,B)=(\frac{n}{d},\frac{\rho}{d})$ if
$d\neq0$, $0\le A_\epsilon,B_\epsilon, A, B\le 1$, and
$(A_\epsilon d_\epsilon,B_\epsilon
d_\epsilon)=(n_\epsilon,\rho_\epsilon), (Ad,Bd)=(n,\rho)$,
$(\rho,n)$ is the limit of $(\rho_{\epsilon},n_{\epsilon})$ in a
suitable weak topology.

For any $\psi\in C([0,t])$, $\phi\in C(\overline{\Omega})$ where
$\psi,\phi\ge0$, we have \be\label{2-II}\begin{split}
&\int_0^t\psi\int_\Omega\phi
(n_\epsilon^\alpha+\rho_\epsilon^\gamma)(\rho_\epsilon+n_\epsilon)
\,dx\,ds\\=& \int_0^t\psi\int_\Omega\phi (A^\alpha
d_\epsilon^\alpha+B^\gamma d_\epsilon^\gamma)(A+B)
d_\epsilon\,dx\,ds\\& +\int_0^t\psi\int_\Omega\phi (A^\alpha
d_\epsilon^\alpha+B^\gamma
d_\epsilon^\gamma)(A_\epsilon-A+B_\epsilon-B)
d_\epsilon\,dx\,ds\\&
 +\int_0^t\psi\int_\Omega\phi \Big[(A_\epsilon^\alpha-A^\alpha)
d_\epsilon^\alpha+(B_\epsilon^\gamma-B^\gamma)
d_\epsilon^\gamma\Big](n_\epsilon+\rho_\epsilon) \,dx\,ds\\=&
\sum\limits_{i=1}^3II_i.
\end{split}
\ee

 For $II_2$, there exists a positive integer $k_0$ large
enough such that \be\label{2-18}
\max\{\frac{k_0\gamma}{k_0-1},\frac{k_0\alpha}{k_0-1}\}\le\beta
\ee due to the assumption that $\max\{\alpha,\gamma\}<\beta$.

Using the H\"older inequality, Lemma \ref{le:h-inofrho} and
(\ref{2-18}), we have \be\label{2-II2}\begin{split} |II_2|
\le&C\left(\int_0^T\int_\Omega
d_\epsilon|A_\epsilon-A|^{k_0}\,dx\,dt\right)^\frac{1}{k_0}\left(\int_0^T\int_\Omega
d_\epsilon|A^\alpha d_\epsilon^\alpha+B^\gamma
d_\epsilon^\gamma|^\frac{k_0}{k_0-1}\,dx\,dt\right)^\frac{k_0-1}{k_0}\\&+
C\left(\int_0^T\int_\Omega
d_\epsilon|B_\epsilon-B|^{k_0}\,dx\,dt\right)^\frac{1}{k_0}\left(\int_0^T\int_\Omega
d_\epsilon|A^\alpha d_\epsilon^\alpha+B^\gamma
d_\epsilon^\gamma|^\frac{k_0}{k_0-1}\,dx\,dt\right)^\frac{k_0-1}{k_0}\\
\le&C\left(\int_0^T\int_\Omega
d_\epsilon|A_\epsilon-A|^{k_0}\,dx\,dt\right)^\frac{1}{k_0}+
C\left(\int_0^T\int_\Omega
d_\epsilon|B_\epsilon-B|^{k_0}\,dx\,dt\right)^\frac{1}{k_0}.
\end{split}
\ee
Choosing $\nu_k=\epsilon$ for Theorem \ref{main 2-le:important},
we conclude that
 \be\label{4.14}\begin{split}
 &\left(\int_0^T\int_\Omega d_\epsilon|A_\epsilon-A|^{k_0}\,dx\,dt\right)^\frac{1}{k_0}\to
 0,\\
 &\left(\int_0^T\int_\Omega d_\epsilon|B_\epsilon-B|^{k_0}\,dx\,dt\right)^\frac{1}{k_0}\to 0
 \end{split}
 \ee
as $\epsilon$ goes to zero. In fact, $d_{\epsilon}\in
L^{\infty}(0,T;L^{\beta}(\Omega))$
 for $\beta>4$, and $u_{\epsilon}\in L^2(0,T;H^1_0(\Omega)),$ and
 \begin{equation*}
  \sqrt{\epsilon}\|\nabla\rho_{\epsilon}\|_{L^2(0,T;L^2(\Omega))}\leq C_0,\,\sqrt{\epsilon}\|\nabla n_{\epsilon}\|_{L^2(0,T;L^2(\Omega))}
  \leq C_0,
\end{equation*}
 and for any $\epsilon>0$ and any $t> 0$:
\begin{equation}
\label{initial condition for n2 over density}\int_{\Omega}
\frac{b_{\epsilon}^2}{d_{\epsilon}}\,dx\leq
\int_{\Omega}\frac{b_0^2}{d_0}\,dx
\end{equation} where $d_\epsilon=\rho_\epsilon+n_\epsilon$, $b_\epsilon=\rho_\epsilon$,
$n_\epsilon$, and (\ref{initial condition for n2 over density}) is
obtained in Remark \ref{remark for an inequality}. Thus, we are
able to apply Theorem \ref{main 2-le:important} to deduce
\eqref{4.14}. Hence we have $II_2\to 0$ as $\epsilon\to 0$.

 For $II_3$, there exists a positive integer $k_1$ large
enough such that \be\label{2-au}\begin{split}
(\alpha+1-\frac{1}{k_1})\frac{k_1}{k_1-1}<\beta+1,\\
(\gamma+1-\frac{1}{k_1})\frac{k_1}{k_1-1}<\beta+1,
\end{split}\ee due to the assumption $\alpha<\beta$. We employ the
H\"older inequality to have \be\label{2-12}
\begin{split}
 |II_3|\le&
C\left(\int_0^T\int_\Omega
d_\epsilon^{(\alpha+1-\frac{1}{k_1})\frac{k_1}{k_1-1}}\,dx\,dt\right)^\frac{k_1-1}{k_1}
\left(\int_0^T\int_\Omega d_\epsilon
\big|A_\epsilon^\alpha-A^\alpha\big|^{k_1}\,dx\,dt\right)^\frac{1}{k_1}\\&+
C\left(\int_0^T\int_\Omega
d_\epsilon^{(\gamma+1-\frac{1}{k_1})\frac{k_1}{k_1-1}}\,dx\,dt\right)^\frac{k_1-1}{k_1}
\left(\int_0^T\int_\Omega d_\epsilon
\big|B_\epsilon^\gamma-B^\gamma\big|^{k_1}\,dx\,dt\right)^\frac{1}{k_1}
\\ \le& C\left(\int_0^T\int_\Omega d_\epsilon
\big|A_\epsilon^\alpha-A^\alpha\big|^{k_1}\,dx\,dt\right)^\frac{1}{k_1}+
C\left(\int_0^T\int_\Omega d_\epsilon
\big|B_\epsilon^\gamma-B^\gamma\big|^{k_1}\,dx\,dt\right)^\frac{1}{k_1}
\rightarrow 0
\end{split}
\ee as $\epsilon\rightarrow0^+$, where we have used (\ref{nrho}),
 (\ref{2-au}), Lemma \ref{le:h-inofrho}, and the fact that \begin{equation}
\begin{split}
\label{control by theorem} \int_0^T\int_\Omega d_\epsilon
\big|A_\epsilon^\alpha-A^\alpha\big|^{k_1}\,dx\,dt
 \le&\alpha^{k_1}\int_0^T\int_\Omega d_\epsilon
\big(\max\{A_\epsilon,A\}\big)^{\alpha-1}\big|A_\epsilon-A\big|^{k_1}\,dx\,dt
\\ \le&C\int_0^T\int_\Omega d_\epsilon
\big|A_\epsilon-A\big|^{k_1}\,dx\,dt\rightarrow 0,\\[2mm]
\int_0^T\int_\Omega d_\epsilon
\big|B_\epsilon^\gamma-B^\gamma\big|^{k_1}\,dx\,dt
 \le&\gamma^{k_1}\int_0^T\int_\Omega d_\epsilon
\big(\max\{B_\epsilon,B\}\big)^{\gamma-1}\big|B_\epsilon-B\big|^{k_1}\,dx\,dt
\\ \le&C\int_0^T\int_\Omega d_\epsilon
\big|B_\epsilon-B\big|^{k_1}\,dx\,dt\rightarrow 0
\end{split}
\end{equation}
 as $\epsilon\rightarrow0^+$, due to  Theorem \ref{main 2-le:important}
 with $\nu_K=\epsilon.$

By virtue of (\ref{2-II}), (\ref{2-II2}) and (\ref{2-12}), one
deduces that \be\label{2-II+1}\begin{split}
\lim\limits_{\epsilon\rightarrow0^+}\int_0^t\psi\int_\Omega\phi
(n_\epsilon^\alpha+\rho_\epsilon^\gamma)(\rho_\epsilon+n_\epsilon)
\,dx\,ds=&
\int_0^t\psi\int_\Omega\phi(A+B)(A^\alpha\overline{d^{\alpha+1}}+B^\gamma\overline{d^{\gamma+1}})\,dx\,ds\\
\ge& \int_0^t\psi\int_\Omega\phi(A+B)(A^\alpha \overline{d^{\alpha}}d+B^\gamma \overline{d^{\gamma}}d)\,dx\,ds\\
=& \int_0^t\psi\int_\Omega\phi\, d(A^\alpha
\overline{d^{\alpha}}+B^\gamma \overline{d^{\gamma}})\,dx\,ds
\end{split}
\ee  where we have used $A+B=1$, $\overline{d^{\alpha+1}}\ge
\overline{d^\alpha}d$, and
$\overline{d^{\gamma+1}}\ge\overline{d^\gamma}d$, because the
functions $z\mapsto z^\alpha (\mathrm{or}\, z^\gamma)$ and
$z\mapsto z$ are increasing functions.

On the other hand, \be\label{2-19}
\begin{split}
&\int_0^t\psi\int_\Omega\phi
(\rho+n)\overline{n^\alpha+\rho^\gamma}\,dx\,ds
\\=&\lim\limits_{\epsilon\rightarrow0^+}\int_0^t\psi\int_\Omega\phi
(\rho+n) (n_\epsilon^\alpha+\rho_\epsilon^\gamma)\,dx\,ds
\\
=&\lim\limits_{\epsilon\rightarrow0^+}\int_0^t\psi\int_\Omega\phi
d (A^\alpha d_\epsilon^\alpha+B^\gamma
d_\epsilon^\gamma)\,dx\,ds+\lim\limits_{\epsilon\rightarrow0^+}\int_0^t\psi\int_\Omega\phi
d
d_\epsilon^\alpha(A_\epsilon^\alpha-A^\alpha)\,dx\,ds\\&+\lim\limits_{\epsilon\rightarrow0^+}\int_0^t\psi\int_\Omega\phi
d d_\epsilon^\gamma(B_\epsilon^\gamma-B^\gamma)\,dx\,ds
\\
=&\int_0^t\psi\int_\Omega\phi d(A^\alpha
\overline{d^{\alpha}}+B^\gamma \overline{d^{\gamma}})\,dx\,ds,
\end{split}
\ee thanks to
\bex\begin{split}\lim\limits_{\epsilon\rightarrow0^+}\int_0^t\psi\int_\Omega\phi\,d
d_\epsilon^\alpha(A_\epsilon^\alpha-A^\alpha)\,dx\,ds\to 0,\quad
\mathrm{and}\quad
\lim\limits_{\epsilon\rightarrow0^+}\int_0^t\psi\int_\Omega\phi\,d
d_\epsilon^\gamma(B_\epsilon^\gamma-B^\gamma)\,dx\,ds\to 0,
\end{split}
\eex
 as
$\epsilon\rightarrow0^+$.

By (\ref{2-II+1}) and (\ref{2-19}), we complete the proof of the
lemma.

\endpf

Since $\psi$ and $\phi$ are arbitrary, we immediately get
\begin{corollary}
Let $(\rho_\epsilon,n_\epsilon)$ be the solution stated in Proposition
\ref{0-le:aweak solution}, and $(\rho, n)$ be the limit in the
sense of (\ref{2-lim}), then \be\label{2-a1}  (\rho+n)
\overline{n^\alpha+\rho^\gamma} \le \overline{(\rho+n)
(n^\alpha+\rho^\gamma)}\ee a.e. on $\Omega\times(0,T)$.
\end{corollary}

Now we are ready to control the right hand side of (\ref{2-17}) in the following lemma.

\begin{lemma}Let $(\rho_\epsilon,n_\epsilon)$ be the solution stated in Lemma
\ref{0-le:aweak solution}, and $(\rho, n)$ be the limit in the
sense of (\ref{2-lim}), then \be\label{2-a2}\int_0^t\int_\Omega
(\rho+n)
\mathrm{div}u\,dx\,ds\le\lim\limits_{\epsilon\rightarrow0^+}\int_0^t\int_\Omega
(\rho_\epsilon+n_\epsilon) \mathrm{div}u_\epsilon\,dx\,ds \ee for
a.e. $t\in(0,T)$.
\end{lemma}
\pf For $\psi_j\in C_0^\infty(0,t)$, $\phi_j\in
C_0^\infty(\Omega)$ given by \be\label{2-test1}\psi_j\in
C_0^\infty(0,T),\quad \psi_j(t)\equiv 1\ \mathrm{for}\
\mathrm{any}\ t\in[\frac{1}{j},T-\frac{1}{j}],\ 0\le\psi_j\le1,\
\psi_j\rightarrow 1,\ee as $j\rightarrow\infty$, and
\be\label{2-test2}\phi_j\in C_0^\infty(\Omega),\quad
\phi_j(x)\equiv1\ \mathrm{for}\ \mathrm{any}\ x\in\big\{x\in\Omega
\big| \mathrm{dist(x,\partial\Omega)\ge\frac{1}{j}}\big\},\
0\le\phi_j\le1,\ \phi_j\rightarrow 1,\ee as $j\rightarrow\infty$,
respectively, then
\be\label{2-RHS}\begin{split}&\int_0^t\int_\Omega (\rho+n)
\mathrm{div}u\,dx\,ds\\=&\int_0^t\psi_j\int_\Omega\phi_j (\rho+n)
\mathrm{div}u\,dx\,ds+\int_0^t\int_\Omega(1-\psi_j\phi_j) (\rho+n)
\mathrm{div}u\,dx\,ds\\=&\frac{1}{2\mu+\lambda}\int_0^t\psi_j\int_\Omega\phi_j
(\rho+n)
\overline{n^\alpha+\rho^\gamma}\,dx\,ds+\frac{1}{2\mu+\lambda}\int_0^t\psi_j\int_\Omega\phi_j
(\rho+n)
\overline{\delta(\rho+n)^\beta}\,dx\,ds\\&-\frac{1}{2\mu+\lambda}\int_0^t\psi_j\int_\Omega\phi_j
(\rho+n) H\,dx\,ds+\int_0^t\int_\Omega(1-\psi_j\phi_j) (\rho+n)
\mathrm{div}u\,dx\,ds\\=& RHS_1+RHS_2+RHS_3+RHS_4,
\end{split}
\ee where we have used
$$(2\mu+\lambda)\mathrm{div}u=\overline{n^\alpha+\rho^\gamma+\delta(\rho+n)^\beta}-H,$$
and
$$
\overline{n^\alpha+\rho^\gamma+\delta(\rho+n)^\beta}=\overline{n^\alpha+\rho^\gamma}+\overline{\delta(\rho+n)^\beta}.$$

For $RHS_2$, we have \be\label{2-RHS2}
\begin{split}
RHS_2=&\frac{1}{2\mu+\lambda}\int_0^t\psi_j\int_\Omega\phi_j
(\rho+n) \overline{\delta(\rho+n)^\beta}\,dx\,ds\\ \le&
\frac{1}{2\mu+\lambda}\liminf\limits_{\epsilon\rightarrow0^+}\int_0^t\psi_j\int_\Omega\phi_j
\delta(\rho_\epsilon+n_\epsilon)
(\rho_\epsilon+n_\epsilon)^\beta\,dx\,ds,
\end{split}
\ee becuase $z\mapsto z$ and $z\mapsto z^\beta$ are
increasing functions.

By virtue of (\ref{2-RHS}), (\ref{2-a1}), (\ref{2-RHS2}), and
(\ref{2-8}), we have
\begin{equation}\begin{split}\label{sss}&\int_0^t\int_\Omega
(\rho+n) \mathrm{div}u\,dx\,ds\\
\le&\frac{1}{2\mu+\lambda}\lim\limits_{\epsilon\rightarrow0^+}\int_0^t\psi_j\int_\Omega\phi_j
(n_\epsilon^\alpha+\rho_\epsilon^\gamma)(\rho_\epsilon+n_\epsilon)
\,dx\,ds\\&+\frac{1}{2\mu+\lambda}\liminf\limits_{\epsilon\rightarrow0^+}\int_0^t\psi_j\int_\Omega\phi_j
(\rho_\epsilon+n_\epsilon)
\delta(\rho_\epsilon+n_\epsilon)^\beta\,dx\,ds\\&-\frac{1}{2\mu+\lambda}\lim\limits_{\epsilon\rightarrow0^+}\int_0^t\psi_j\int_\Omega\phi_j
(\rho_\epsilon+n_\epsilon)
H_\epsilon\,dx\,ds+\int_0^t\int_\Omega(1-\psi_j\phi_j) (\rho+n)
\mathrm{div}u\,dx\,ds\\
\le&\lim\limits_{\epsilon\rightarrow0^+}\int_0^t\int_\Omega
(\rho_\epsilon+n_\epsilon)
\mathrm{div}u_\epsilon\,dx\,ds+\lim\limits_{\epsilon\rightarrow0^+}\int_0^t\int_\Omega(\psi_j\phi_j-1)
(\rho_\epsilon+n_\epsilon)
\mathrm{div}u_\epsilon\,dx\,ds\\&+\int_0^t\int_\Omega(1-\psi_j\phi_j)
(\rho+n) \mathrm{div}u\,dx\,ds.
\end{split}
\end{equation} Letting $j\rightarrow\infty$ in \eqref{sss}, we complete the proof of the
lemma.
\endpf
\begin{flushleft}
\textbf{    Step 3: Strong convergence of $\rho_{\epsilon}$ and $n_{\epsilon}$}
\end{flushleft}
The  main task is to show the strong convergence of $\rho_{\epsilon}$ and $n_{\epsilon}$. This yields Proposition \ref{key pro}. In particular,
 With (\ref{2-a2}), letting
$\epsilon\rightarrow0^+$ in \eqref{2-17}, we deduce that \bex
\begin{split}
\int_\Omega \big[\overline{\rho\log \rho}-\rho\log
\rho+\overline{n\log n}-n\log n\big](t)\,dx \le0.
\end{split}
\eex Thanks to the convexity of $z\mapsto z\log z$, we have \bex
\overline{\rho\log \rho}\ge\rho\log \rho\quad \mathrm{and}\quad
\overline{n\log n}\ge n \log n \eex a.e. on $Q_T$. This turns out
that \bex \int_\Omega \big[\overline{\rho\log \rho}-\rho\log
\rho+\overline{n\log n}-n\log n\big](t)\,dx=0. \eex Hence we get
\bex \overline{\rho\log \rho}=\rho\log \rho\quad \mathrm{and}\quad
\overline{n\log n}= n \log n \eex a.e. on $Q_T$, which combined
with Lemma \ref{le:h-inofrho} implies strong convergence of
$\rho_\epsilon,n_\epsilon$ in $L^\beta(Q_T)$. Thus we complete the
proof.

\bigskip

With Proposition \ref{key pro},  we recover a global weak solution
to the system (\ref{a2-equation}), (\ref{a2-initial}) and
(\ref{a2-boundary}) with
$\overline{n^\alpha+\rho^\gamma+\delta(\rho+n)^\beta}$ replaced by
$n^\alpha+\rho^\gamma+\delta(\rho+n)^\beta$.

\begin{proposition}\label{2-le:aweak solution}
Suppose $\beta>\max\{4,\alpha,\gamma\}$. For any given $\delta>0$,
there exists a global weak solution
$(\rho_\delta,n_\delta,u_\delta)$ to the following system over
$\Omega\times(0,\infty)$:
 \begin{equation}\label{a3-equation}
    \left\{
    \begin{array}{l}
        n_t+\mathrm{div}(n u)=0,\\[2mm]
        \rho_t+\mathrm{div}(\rho u)=0,\\[2mm]
 \big[(\rho+n) u\big]_t+\mathrm{div}\big[(\rho+n) u\otimes u\big]
 +\nabla(n^\alpha+\rho^\gamma+\delta(\rho+n)^\beta) =
            \mu \Delta u + (\mu+\lambda)\nabla \mathrm{div}u,
    \end{array}
    \right.
    \end{equation}
    with initial and boundary condition
\be\label{a3-initial}\big (\rho, n,(\rho+n)
u\big)|_{t=0}=(\rho_{0,\delta},n_{0,\delta},M_{0,\delta})\
\mathrm{on}\ \overline{\Omega}, \ee \be\label{a3-boundary}
u|_{\partial\Omega}=0 \ \mathrm{for}\ t\ge0, \ee such that for any
given $T>0$, the following estimates
 \be\label{2-r1}
\sup\limits_{t\in[0,T]}\|\rho_\delta(t)\|_{L^\gamma(\Omega)}^\gamma\le
C(\rho_0,n_0,M_0),\ee

\be\label{2-r1+1}
\sup\limits_{t\in[0,T]}\|n_\delta(t)\|_{L^\alpha(\Omega)}^\alpha\le
C(\rho_0,n_0,M_0),\ee

\be\label{2-r2}\delta\sup\limits_{t\in[0,T]}\|(\rho_\delta(t),n_\delta(t))\|_{L^\beta(\Omega)}^\beta\le
C(\rho_0,n_0,M_0),\ee

\be\label{2-r3}\sup\limits_{t\in[0,T]}\|\sqrt{\rho_\delta+n_\delta}(t)u_\delta(t)\|_{L^2(\Omega)}^2\le
C(\rho_0,n_0,M_0),\ee

\be\label{2-r4} \int_0^T\|u_\delta(t)\|_{H^1_0(\Omega)}^2\,dt\le
C(\rho_0,n_0,M_0),\ee and \be\label{2-r5}
\|(\rho_\delta(t),n_\delta(t))\|_{L^{\beta+1}(Q_T)}\le
C(\beta,\delta,\rho_0,n_0,M_0) \ee hold, where the norm
$\|(\cdot,\cdot)\|$ denotes $\|\cdot\| + \|\cdot\|$. Besides, if
$\frac{1}{c_0}\rho_0\le n_0\le c_0\rho_0$, we have
\be\label{3-nrho}\begin{cases} \frac{1}{c_0}\rho_\delta(x,t)\le
n_\delta(x,t)\le c_0\rho_\delta(x,t)\quad \mathrm{for}\
\mathrm{a.e.}\ (x,t)\in Q_T,
\\[2mm]
\sup\limits_{t\in[0,T]}\|(\rho_\delta,n_\delta)(t)\|_{L^{\alpha_1}(\Omega)}^{\alpha_1}\le
C(\rho_0,n_0,M_0),\end{cases}\ee where
$\alpha_1=\max\{\alpha,\gamma\}$.
\end{proposition}

\section{Passing to the limit in the artificial pressure term as
$\delta\rightarrow0^+$} \setcounter{equation}{0}
\setcounter{theorem}{0}

In this section,  we shall recover the weak solution to \eqref{equation}-\eqref{boundary} by passing to the limit of
$(\rho_\delta,n_\delta,u_\delta)$ as
$\delta\rightarrow0$. Note
that the estimate (\ref{2-r5}) depends on $\delta$. Thus to begin
with, we have to get the higher integrability estimates of the
pressure term uniformly for $\delta$. Let $C$ be a generic constant independent of $\delta$ which will
be used throughout this section.

\subsection{Passing to the limit as $\delta\rightarrow0^+$}

We can follow the similar argument as in \cite{Feireisl} to have
the higher integrability estimates of $\rho$ and $n$ in the
following lemma. We only need to
 modify the proof a little bit on $n$.

\begin{lemma}\label{3-le:h-inofrho} Let $(\rho_\delta,n_\delta,u_\delta)$ be the solution
stated in Proposition \ref{2-le:aweak solution}, then
\be\label{3-hierho} \int_0^T\int_\Omega
(n_\delta^{\alpha+\theta_1}+\rho_\delta^{\gamma+\theta_2}+\delta
n_\delta^{\beta+\theta_1}+\delta\rho_\delta^{\beta+\theta_2})\,dx
dt\le C(\theta_1,\theta_2)\ee for any positive constants
$\theta_1$ and $\theta_2$ satisfying
\bex\begin{cases}\theta_1<\frac{\alpha}{3}\,\, \mathrm{and}\,\,
\theta_1\le\min\{1,\frac{2\alpha}{3}-1\};\,\,
\theta_2<\frac{\gamma}{3}\,\, \mathrm{and}\,\,
\theta_2\le\min\{1,\frac{2\gamma}{3}-1\}\quad \mathrm{if}\,\,
\alpha,\gamma\in (\frac{3}{2},\infty),\\[2mm]
\theta<\frac{\max\{\alpha,\gamma\}}{3}\, \,\mathrm{and}\,\,
\theta\le\min\{1,\frac{2\max\{\alpha,\gamma\}}{3}-1\}\quad
\mathrm{if}\,\,\gamma\in (\frac{3}{2},\infty),\,\,
\alpha\in[1,\infty),\,\, \mathrm{and}\,\, \frac{1}{c_0}\rho_0\le
n_0\le c_0\rho_0,
\end{cases}\eex where $\theta=\theta_1=\theta_2$.
\end{lemma}

In view of (\ref{3-hierho}) and (\ref{3-nrho}), we have the
following corollary.

\begin{corollary}
Let $(\rho_\delta,n_\delta,u_\delta)$ be the solution stated in
Proposition \ref{2-le:aweak solution}, if $\frac{1}{c_0}\rho_0\le
n_0\le c_0\rho_0$, then \be\label{2-r6}\int_0^T\int_\Omega
\rho_\delta^{\max\{\alpha+\theta_1,\gamma+\theta_2\}}+n_\delta^{\max\{\alpha+\theta_1,\gamma+\theta_2\}}\,dx\,dt
\le C.\ee
\end{corollary}

With (\ref{2-r1}), (\ref{2-r1+1}), (\ref{2-r3}), (\ref{2-r4}),
(\ref{3-nrho}), (\ref{3-hierho}), and (\ref{2-r6}), letting
$\delta\rightarrow0^+$ (take the subsequence if necessary), we
have

{\noindent\bf Case 1.} $\alpha,\gamma\in(\frac{9}{5},\infty)$, and
$\max\{\frac{3\gamma}{4},\gamma-1,\frac{3(\gamma+1)}{5}\}<\alpha<\min\{\frac{4\gamma}{3},\gamma+1,\frac{5\gamma}{3}-1\}$.
 \be\label{3-lim}
\begin{cases}
\rho_\delta\rightarrow \rho\ \mathrm{in}\ C([0,T];
L_{weak}^\gamma(\Omega))\ \mathrm{and}\ \mathrm{weakly}\
\mathrm{in}\ L^{\gamma+\theta_2}(Q_T)\ \mathrm{as}\
\delta\rightarrow0^+,\\[2mm]
n_\delta\rightarrow n\ \mathrm{in}\ C([0,T];
L_{weak}^\alpha(\Omega))\ \ \mathrm{and}\ \mathrm{weakly}\
\mathrm{in}\ L^{\alpha+\theta_1}(Q_T)\ \mathrm{as}\
\delta\rightarrow0^+,\\[2mm]
u_\delta\rightarrow u\ \mathrm{weakly}\ \mathrm{in}\ L^2(0,T;
H_0^1(\Omega))\ \mathrm{as}\ \delta\rightarrow0^+,\\[2mm]
(\rho_\delta+n_\delta) u_\delta\rightarrow (\rho+n) u \
\mathrm{in}\ C([0,T];
L_{weak}^\frac{2\min\{\gamma,\alpha\}}{\min\{\gamma,\alpha\}+1})\cap
C([0,T];H^{-1}(\Omega))\ \mathrm{as}\ \delta\rightarrow0^+,\\[2mm]
(\rho_\delta u_\delta, n_\delta u_\delta)\rightarrow (\rho u, n u)
\ \mathrm{in}\ \mathcal{D}^\prime(Q_T)\ \mathrm{as}\
\delta\rightarrow0^+,\\[2mm]
(\rho_\delta+n_\delta) u_\delta\otimes u_\delta\rightarrow
(\rho+n) u\otimes u\ \mathrm{in}\ \mathcal{D}^\prime(Q_T)\
\mathrm{as}\
\delta\rightarrow0^+,\\[2mm]
n_\delta^\alpha+\rho_\delta^\gamma\rightarrow
\overline{n^\alpha+\rho^\gamma}\ \mathrm{weakly}\ \mathrm{in}\
L^{\min\{\frac{\gamma+\theta_2}{\alpha},\frac{\alpha+\theta_1}{\gamma}\}}(Q_T)\
\mathrm{as}\
\delta\rightarrow0^+,\\[2mm]
\delta(\rho_\delta+n_\delta)^\beta\rightarrow 0\ \mathrm{in}\
L^1(Q_T)\ \mathrm{as}\ \delta\rightarrow0^+.
\end{cases}
\ee {\noindent\bf Case 2.} $\alpha\in[1,\infty),\,\,
\gamma\in(\frac{9}{5},\infty)$, and $\frac{1}{c_0}\rho_0\le n_0\le
c_0\rho_0$.

\be\label{3-lim1}
\begin{cases}
(\rho_\delta,n_\delta)\rightarrow (\rho,n)\ \mathrm{in}\ C([0,T];
L_{weak}^{\max\{\gamma,\alpha\}}(\Omega))\ \mathrm{and}\
\mathrm{weakly}\ \mathrm{in}\
L^{\max\{\alpha+\theta_1,\gamma+\theta_2\}}(Q_T)\ \mathrm{as}\
\delta\rightarrow0^+,\\[2mm]
u_\delta\rightarrow u\ \mathrm{weakly}\ \mathrm{in}\ L^2(0,T;
H_0^1(\Omega))\ \mathrm{as}\ \delta\rightarrow0^+,\\[2mm]
(\rho_\delta+n_\delta) u_\delta\rightarrow (\rho+n) u \
\mathrm{in}\ C([0,T];
L_{weak}^\frac{2\max\{\gamma,\alpha\}}{\max\{\gamma,\alpha\}+1})\cap
C([0,T];H^{-1}(\Omega))\ \mathrm{as}\ \delta\rightarrow0^+,\\[2mm]
(\rho_\delta u_\delta, n_\delta u_\delta)\rightarrow (\rho u, n u)
\ \mathrm{in}\ \mathcal{D}^\prime(Q_T)\ \mathrm{as}\
\delta\rightarrow0^+,\\[2mm]
(\rho_\delta+n_\delta) u_\delta\otimes u_\delta\rightarrow
(\rho+n) u\otimes u\ \mathrm{in}\ \mathcal{D}^\prime(Q_T)\
\mathrm{as}\
\delta\rightarrow0^+,\\[2mm]
n_\delta^\alpha+\rho_\delta^\gamma\rightarrow
\overline{n^\alpha+\rho^\gamma}\ \mathrm{weakly}\ \mathrm{in}\
L^{\max\{\frac{\gamma+\theta_2}{\alpha},\frac{\alpha+\theta_1}{\gamma}\}}(Q_T)\
\mathrm{as}\
\delta\rightarrow0^+,\\[2mm]
\delta(\rho_\delta+n_\delta)^\beta\rightarrow 0\ \mathrm{in}\
L^1(Q_T)\ \mathrm{as}\ \delta\rightarrow0^+,\\[2mm]
\frac{1}{c_0}\rho_\delta(x,t)\le n_\delta(x,t)\le
c_0\rho_\delta(x,t),\,\, \mathrm{for}\,\, \mathrm{a.e.,}\,\,
(x,t)\in Q_T,\\[2mm]
\frac{1}{c_0}\rho(x,t)\le n(x,t)\le c_0\rho(x,t),\,\,
\mathrm{for}\,\, \mathrm{a.e.,}\,\, (x,t)\in Q_T.
\end{cases}
\ee

In summary, the limit $(\rho, n, u)$ solves the following system
in the sense of distribution over $\Omega\times [0, T]$ for any
given $T>0$:
 \begin{equation}\label{3-equation}
    \left\{
    \begin{array}{l}
        n_t+\mathrm{div}(n u)=0,\\[2mm]
        \rho_t+\mathrm{div}(\rho u)=0,\\[2mm]
 \big[(\rho+n) u\big]_t+\mathrm{div}\big[(\rho+n) u\otimes u\big]
 +\nabla(\overline{\rho^\gamma+n^\alpha}) =
            \mu \Delta u + (\mu+\lambda)\nabla \mathrm{div}u,
    \end{array}
    \right.
    \end{equation}
    with initial and boundary condition
\be\label{3-initial} (\rho, n,(\rho+n)
u)|_{t=0}=(\rho_0,n_0,M_0)\quad \mathrm{on}\ \overline{\Omega},
\ee \be\label{3-boundary} u|_{\partial\Omega}=0 \quad
\mathrm{for}\ t\ge0, \ee where the convergence of the approximate
initial data in (\ref{a3-initial}) is due to (\ref{appinitial}).

To recover a weak solution to \eqref{equation}-\eqref{boundary}, we only need to show the following claim:
\begin{itemize}
\item {\bf\large Claim.
$\overline{\rho^\gamma+n^\alpha}=\rho^\gamma+n^\alpha$.}
\end{itemize}

\subsection{The weak limit of pressure}
The objective of this subsection is to show the strong convergence
of $\rho_{\delta}$ and $n_{\delta}$ as $\delta$ goes to zero. This
allows us to  prove
$\overline{\rho^\gamma+n^\alpha}=\rho^\gamma+n^\alpha$ as
$\delta\to 0$. From now, we need that $\rho_{\delta}$ is bounded
in $L^q(Q_T)$ for some $q>2$. By Lemma
\ref{3-le:h-inofrho}, we need the restriction
$\gamma>\frac{9}{5}.$ 

We consider a family of cut-off functions introduced in
\cite{Feireisl} and references therein, i.e., \bex
T_k(z)=kT(\frac{z}{k}),\ z\in\mathbb{R},\ k=1,2,\cdot\cdot\cdot
\eex where $T\in C^\infty(\mathbb{R})$ satisfying \bex T(z)=z\
for\ z\le1,\quad T(z)=2\ for\ z\ge3.\quad T\ is\ concave. \eex The
cut-off functions will be used in particular to handle the cross
terms due to the two-variable pressure, see the proof of Lemma
\ref{3-le:4.6}. Since $\rho_\delta\in L^2(Q_T),u_\delta\in
L^2(0,T; H_0^1(\Omega))$, Lemma \ref{main lemma} suggests that
$(\rho_\delta,u_\delta)$ is a renormalized solution of
(\ref{3-equation})$_2$. Thus we have \be\label{3-Tk}
[T_k(\rho_\delta)]_t+\mathrm{div}[T_k(\rho_\delta)u_\delta]
+[T_k^\prime(\rho_\delta)\rho_\delta-T_k(\rho_\delta)]\mathrm{div}u_\delta=0\quad
\mathrm{in}\ \mathcal{D}^\prime(Q_T). \ee

For any given $k$, $T_k(\rho_\delta)$ is bounded in
$L^\infty(Q_T)$. Passing to the limit as $\delta\rightarrow0^+$
(taking the subsequence if necessary), we have \bex\begin{split}
&T_k(\rho_\delta)\rightarrow\overline{T_k(\rho)}\
\mathrm{weak*}\ \mathrm{in}\ L^\infty(Q_T),\\
&T_k(\rho_\delta)\rightarrow \overline{T_k(\rho)}\ \mathrm{in}\
C([0,T]; L^p_{weak}(\Omega)),\ \mathrm{for}\ \mathrm{any}\
p\in[1,\infty),
\\ &T_k(\rho_\delta)\rightarrow \overline{T_k(\rho)}\ \mathrm{in}\
C([0,T];H^{-1}(\Omega)),\\
&[T_k^\prime(\rho_\delta)\rho_\delta-T_k(\rho_\delta)]\mathrm{div}u_\delta\rightarrow
\overline{[T_k^\prime(\rho)\rho-T_k(\rho)]\mathrm{div}u}\
\mathrm{weakly}\ \mathrm{in}\ L^2(Q_T).
\end{split}
\eex This yields \be\label{3-Trho}
[\overline{T_k(\rho)}]_t+\mathrm{div}[\overline{T_k(\rho)}u]+\overline{[T_k^\prime(\rho)\rho-T_k(\rho)]\mathrm{div}u}=0\
\mathrm{in}\ \mathcal{D}^\prime(\mathbb{R}^3\times(0,T)). \ee
Similarly, we have \be\label{3-Tk-t}
[T_k(n_\delta)]_t+\mathrm{div}[T_k(n_\delta)u_\delta]
+[T_k^\prime(n_\delta)n_\delta-T_k(n_\delta)]\mathrm{div}u_\delta=0\quad
\mathrm{in}\ \mathcal{D}^\prime(Q_T), \ee and
 \be\label{3-Trho-t}
[\overline{T_k(n)}]_t+\mathrm{div}[\overline{T_k(n)}u]+\overline{[T_k^\prime(n)n-T_k(n)]\mathrm{div}u}=0\
\mathrm{in}\ \mathcal{D}^\prime(\mathbb{R}^3\times(0,T)).\ee

 Denote \bex\begin{split}
H_\delta:=&\rho_\delta^\gamma+n_\delta^\alpha-(2\mu+\lambda)\mathrm{div}u_\delta,
\\
\overline{H}:=&\overline{\rho^\gamma+n^\alpha}-(2\mu+\lambda)\mathrm{div}u.
\end{split}
\eex
We will have the following Lemma on $H_{\delta}$ and $\overline{H}$.
\begin{lemma}\label{3-le:3.7} Let
$(\rho_\delta,n_\delta,u_\delta)$ be the solution stated in Proposition
\ref{2-le:aweak solution} and $(\rho,n,u)$ be the limit, then
\be\label{3-8}
\lim\limits_{\delta\rightarrow0^+}\int_0^T\psi\int_\Omega\phi
H_\delta
\big[T_k(\rho_\delta)+T_k(n_\delta)\big]\,dx\,dt=\int_0^T\psi\int_\Omega\phi
\overline{H}\ \big[\overline{T_k(\rho)}+
\overline{T_k(n)}\big]\,dx\,dt, \ee for any $\psi\in
C_0^\infty(0,T)$ and $\phi\in
 C_0^\infty(\Omega)$.
\end{lemma}
\pf The proof is similar to the work of \cite{Feireisl}.

\endpf

Now let us focus on  the following Proposition.
\begin{proposition}
\label{key pro for delta}For any $\alpha,\gamma>\frac{9}{5}$ and
$\max\{\frac{3\gamma}{4},\gamma-1,\frac{3(\gamma+1)}{5}\}<\alpha<\min\{\frac{4\gamma}{3},\gamma+1,\frac{5\gamma}{3}-1\}$,
then
\begin{equation}
\label{3-laststep}
\overline{n^\alpha+\rho^\gamma}=n^\alpha+\rho^\gamma
\end{equation}
a.e. on $Q_T$. In addition, if the initial data satisfy $$\frac{1}{c_0}\rho_0\le n_0\le
c_0\rho_0,$$ then for  $\alpha\geq 1,\,\,
\gamma>\frac{9}{5},$  \eqref{3-laststep} holds.
\end{proposition}
There are two steps to prove it.
\begin{flushleft}
    \textbf{Step 1: Study for the weak limit of $\rho_{\delta}^{\gamma}+n_{\delta}^{\alpha}$}
\end{flushleft}

Relying on Theorem \ref{main 2-le:important} with $\nu_K=0$, we are able to show the following lemma. It is crucial to obtain Proposition \ref{key pro for delta}.

\begin{lemma}\label{3-le:4.6}
Let $(\rho_\delta,n_\delta)$ be the solutions constructed in
Propsotion \ref{2-le:aweak solution}, and $(\rho,n)$ be the limit,
then  \be\label{3-c1}\begin{split} \int_0^t\psi\int_\Omega\phi
\Big[\overline{T_k(\rho)}
+\overline{T_k(n)}\Big]\big(\overline{\rho^\gamma+n^\alpha}\big)\,dx\,ds
 \le\int_0^t\psi\int_\Omega\phi
\overline{\big[T_k(\rho)+T_k(n)\big]\big(\rho^\gamma+n^\alpha\big)}\,dx\,ds.
\end{split}
\ee for any $t\in[0,T]$ and any $\psi\in C[0,t]$, $\phi\in
C(\overline{\Omega})$ where $\psi,\phi\ge0$,
$\alpha,\gamma>\frac{9}{5}$ and
$$\max\{\frac{3\gamma}{4},\gamma-1,\frac{3(\gamma+1)}{5}\}<\alpha<\min\{\frac{4\gamma}{3},\gamma+1,\frac{5\gamma}{3}-1\}.$$
In additional, if initial data satisfies $$\frac{1}{c_0}\rho_0\le n_0\le c_0\rho_0,$$ then for
$\alpha\geq 1,\,\, \gamma>\frac{9}{5},$ \eqref{3-c1} holds.
\end{lemma}

\pf \be\label{3-sum}\begin{split}
&\lim\limits_{\delta\rightarrow0^+}\int_0^t\psi\int_\Omega\phi
\big[T_k(\rho_\delta)+T_k(n_\delta)\big]\big(\rho_\delta^\gamma+n_\delta^\alpha\big)\,dx\,ds\\
=&\lim\limits_{\delta\rightarrow0^+}\int_0^t\psi\int_\Omega\phi
T_k(\rho_\delta)\rho_\delta^\gamma\,dx\,ds+\lim\limits_{\delta\rightarrow0^+}\int_0^t\psi\int_\Omega\phi
T_k(\rho_\delta)n_\delta^\alpha\,dx\,ds\\&+\lim\limits_{\delta\rightarrow0^+}\int_0^t\psi\int_\Omega\phi
T_k(n_\delta)\rho_\delta^\gamma\,dx\,ds+\lim\limits_{\delta\rightarrow0^+}\int_0^t\psi\int_\Omega\phi
T_k(n_\delta)n_\delta^\alpha\,dx\,ds\\=& \sum\limits_{i=1}^4IV_i.
\end{split}
\ee

For $IV_1$, since $z\mapsto T_k(z)$ and  $z\mapsto z^\gamma$ are
increasing functions, we have \be\label{3-c3}
\begin{split}
0\le&\lim\limits_{\delta\rightarrow0^+}\int_0^t\psi\int_\Omega\phi
\big[T_k(\rho_\delta)-T_k(\rho)\big]\big[\rho_\delta^\gamma-\rho^\gamma\big]\,dx\,ds
\\ =& \int_0^t\psi\int_\Omega\phi
\overline{T_k(\rho)\rho^\gamma}\,dx\,ds-
\int_0^t\psi\int_\Omega\phi
\overline{T_k(\rho)}\rho^\gamma\,dx\,ds-\int_0^t\psi\int_\Omega\phi
T_k(\rho)\overline{\rho^\gamma}\,dx\,ds\\&+
\int_0^t\psi\int_\Omega\phi T_k(\rho)\rho^\gamma\,dx\,ds\\ =&
\int_0^t\psi\int_\Omega\phi
\overline{T_k(\rho)\rho^\gamma}\,dx\,ds-\int_0^t\psi\int_\Omega\phi
\overline{T_k(\rho)}\overline{\rho^\gamma}\,dx\,ds\\&+\int_0^t\psi\int_\Omega\phi
\big[\overline{T_k(\rho)}-T_k(\rho)\big]\big(\overline{\rho^\gamma}-\rho^\gamma\big)\,dx\,ds\\
\le& \int_0^t\psi\int_\Omega\phi
\overline{T_k(\rho)\rho^\gamma}\,dx\,ds-\int_0^t\psi\int_\Omega\phi
\overline{T_k(\rho)}\overline{\rho^\gamma}\,dx\,ds
\end{split}
\ee where we have used the fact
$\overline{\rho^\gamma}\ge\rho^\gamma$ and
$\overline{T_k(\rho)}\le T_k(\rho)$, which could be done by the convexity of $z\mapsto z^\gamma$ and the concavity
of $z\mapsto T_k(z)$.

Thanks to (\ref{3-c3}), we have \be\label{3-IV1}
\int_0^t\psi\int_\Omega\phi
\overline{T_k(\rho)}\overline{\rho^\gamma}\,dx\,ds\le
\int_0^t\psi\int_\Omega\phi
\overline{T_k(\rho)\rho^\gamma}\,dx\,ds=\lim\limits_{\delta\rightarrow0^+}\int_0^t\psi\int_\Omega\phi
T_k(\rho_\delta)\rho_\delta^\gamma\,dx\,ds=IV_1.\ee

Similar to (\ref{3-IV1}), we have
 \be\label{3-IV4}
\int_0^t\psi\int_\Omega\phi
\overline{T_k(n)}\overline{n^\alpha}\,dx\,ds\le
\lim\limits_{\delta\rightarrow0^+}\int_0^t\psi\int_\Omega\phi
T_k(n_\delta)n_\delta^\alpha\,dx\,ds=IV_4.\ee

\bigskip

For $IV_2$, we need to discuss the sizes of $\alpha$ and $\gamma$
in order to guarantee the boundedness of $\rho_\delta^\alpha$ and
$n_\delta^\gamma$ in $L^q(Q_T)$ for some $q>1$.

{\bf Case 1.}  $\alpha,\gamma>\frac{9}{5}$ and
$\max\{\frac{3\gamma}{4},\gamma-1,\frac{3(\gamma+1)}{5}\}<\alpha<\min\{\frac{4\gamma}{3},\gamma+1,\frac{5\gamma}{3}-1\}$.

\bigskip

In this case, there exist two positive constants
$\theta_1\in(\gamma-\alpha,\,\min\{\frac{\alpha}{3},1,\frac{2\alpha}{3}-1\})$
and
$\theta_2\in(\alpha-\gamma,\,\min\{\frac{\gamma}{3},1,\frac{2\gamma}{3}-1\})$,
since
$\max\{\frac{3\gamma}{4},\gamma-1,\frac{3(\gamma+1)}{5}\}<\alpha<\min\{\frac{4\gamma}{3},\gamma+1,\frac{5\gamma}{3}-1\}$
implies that
$$
\gamma-\alpha<\min\{\frac{\alpha}{3},1,\frac{2\alpha}{3}-1\},\quad
\mathrm{and}\quad
\alpha-\gamma<\min\{\frac{\gamma}{3},1,\frac{2\gamma}{3}-1\}.
$$ Note that we are able to take $\theta_1$ and $\theta_2$ here
the same as those in Lemma \ref{3-le:h-inofrho}. Then there exists
a positive integer $k_2$ large enough such that
\be\label{3-au}\begin{cases} 0<\frac{\alpha
k_2}{k_2-1}-\frac{1}{k_2-1}<\gamma+\theta_2,\\[2mm]
0<\frac{\alpha k_2}{k_2-1}-\frac{1}{k_2-1}<\alpha+\theta_1.
 \end{cases}\ee In
this case, $d_\delta=\rho_\delta+n_\delta$ is bounded in
$L^{\frac{\alpha k_2}{k_2-1}-\frac{1}{k_2-1}}(Q_T)$. Then
 \be\label{3-IV2}
\begin{split}
IV_2=&\lim\limits_{\delta\rightarrow0^+}\int_0^t\psi\int_\Omega\phi
T_k(\rho_\delta)d_\delta^\alpha(A_\delta^\alpha-A^\alpha)\,dx\,ds
+\lim\limits_{\delta\rightarrow0^+}\int_0^t\psi\int_\Omega\phi
T_k(\rho_\delta)d_\delta^\alpha
A^\alpha\,dx\,ds\\
\ge&-2kC\lim\limits_{\delta\rightarrow0^+}\Big(\int_0^t\int_\Omega
d_\delta\big|A_\delta^\alpha-A^\alpha\big|^{k_2}\,dx\,ds\Big)^\frac{1}{k_2}\Big(\int_0^t\int_\Omega
d_\delta^{\frac{\alpha
k_2}{k_2-1}-\frac{1}{k_2-1}}\,dx\,ds\Big)^\frac{k_2-1}{k_2}\\&
+\lim\limits_{\delta\rightarrow0^+}\int_0^t\psi\int_\Omega\phi
T_k(B d_\delta)d_\delta^\alpha
A^\alpha\,dx\,ds+\lim\limits_{\delta\rightarrow0^+}\int_0^t\psi\int_\Omega\phi
[T_k(B_\delta
d_\delta)-T_k(B d_\delta)]d_\delta^\alpha A^\alpha\,dx\,ds\\
\ge&-2k
C\alpha\lim\limits_{\delta\rightarrow0^+}\Big(\int_0^t\int_\Omega
d_\delta
\Big|\big(\max\{A_\delta,A\}\big)^{\alpha-1}\big|A_\delta-A\big|\Big|^{k_2}\,dx\,dt\Big)^\frac{1}{k_2}
 \\&+\int_0^t\psi\int_\Omega\phi
\overline{T_k(B d)d^\alpha
}A^\alpha\,dx\,ds+\lim\limits_{\delta\rightarrow0^+}\int_0^t\psi\int_\Omega\phi
[T_k(B_\delta
d_\delta)-T_k(B d_\delta)]d_\delta^\alpha A^\alpha\,dx\,ds\\
\ge&-2k
C\alpha\lim\limits_{\delta\rightarrow0^+}\Big(\int_0^t\int_\Omega
d_\delta \big|A_\delta-A\big|^{k_2}\,dx\,dt\Big)^\frac{1}{k_2}
 \\&+\int_0^t\psi\int_\Omega\phi
\overline{T_k(B d)d^\alpha
}A^\alpha\,dx\,ds+\lim\limits_{\delta\rightarrow0^+}\int_0^t\psi\int_\Omega\phi
[T_k(B_\delta d_\delta)-T_k(B d_\delta)]d_\delta^\alpha
A^\alpha\,dx\,ds,
\end{split}
\ee where
$(A_\delta,B_\delta)=(\frac{n_\delta}{d_\delta},\frac{\rho_\delta}{d_\delta})$
with $d_\delta=\rho_\delta+n_\delta$,
$(A,B)=(\frac{n}{d},\frac{\rho}{d})$ with $d=\rho+n$.

 In view of Theorem \ref{main 2-le:important} with $\nu_K=0$, (\ref{3-IV2}), and the arguments
similar to (\ref{3-c3}), we have
 \be\label{3-IV2+1}
\begin{split}
IV_2 \ge&\int_0^t\psi\int_\Omega\phi \overline{T_k(B
d)}\,\,\overline{d^\alpha}A^\alpha\,dx\,ds+\lim\limits_{\delta\rightarrow0^+}\int_0^t\psi\int_\Omega\phi
[T_k(B_\delta d_\delta)-T_k(B d_\delta)]d_\delta^\alpha
A^\alpha\,dx\,ds\\=&\int_0^t\psi\int_\Omega\phi
\overline{T_k(\rho)}\,\,\overline{d^\alpha}A^\alpha\,dx\,ds+\lim\limits_{\delta\rightarrow0^+}\int_0^t\psi\int_\Omega\phi
\Big(T_k(B
d_\delta)-T_k(\rho_\delta)\Big)\,\,\overline{d^\alpha}A^\alpha\,dx\,ds\\&+\lim\limits_{\delta\rightarrow0^+}\int_0^t\psi\int_\Omega\phi
[T_k(\rho_\delta)-T_k(B d_\delta)]d_\delta^\alpha
A^\alpha\,dx\,ds\\=&\int_0^t\psi\int_\Omega\phi
\overline{T_k(\rho)}\,\,\overline{n^\alpha}\,dx\,ds+\lim\limits_{\delta\rightarrow0^+}\int_0^t\psi\int_\Omega\phi
\Big(T_k(B
d_\delta)-T_k(\rho_\delta)\Big)\,\,\overline{d^\alpha}A^\alpha\,dx\,ds+\\&\lim\limits_{\delta\rightarrow0^+}\int_0^t\psi\int_\Omega\phi
\overline{T_k(\rho)}\,\,(d_\delta^\alpha
A^\alpha-n_\delta^\alpha)\,dx\,ds+\lim\limits_{\delta\rightarrow0^+}\int_0^t\psi\int_\Omega\phi
[T_k(\rho_\delta)-T_k(B d_\delta)]d_\delta^\alpha
A^\alpha\,dx\,ds.
\end{split}
\ee

In view of Theorem \ref{main 2-le:important} with $\nu_K=0$, in
particular, of \eqref{key bound}, we have
\be\label{ndconvergence}\begin{cases}
n_\delta-Ad_\delta\rightarrow0 \quad \mathrm{a.e.}\ \mathrm{in}\
Q_T,\\[2mm]\rho_\delta-Bd_\delta\rightarrow0 \quad \mathrm{a.e.}\ \mathrm{in}\
Q_T,
\end{cases}\ee as $\delta\rightarrow0^+$ (take the
subsequence if necessary). (\ref{ndconvergence})$_2$ implies that
\be\label{3-a4}
T_k\big(Bd_\delta\big)-T_k\big(\rho_\delta\big)\rightarrow0 \quad
\mathrm{a.e.}\ \mathrm{in}\ Q_T\ee as $\delta\rightarrow0^+$ (take
the subsequence if necessary).

Since $\Big(T_k(B
d_\delta)-T_k(\rho_\delta)\Big)\,\,\overline{d^\alpha}A^\alpha$,
$\overline{T_k(\rho)}\,\,(d_\delta^\alpha
A^\alpha-n_\delta^\alpha)$, and $[T_k(\rho_\delta)-T_k(B
d_\delta)]d_\delta^\alpha A^\alpha$ are bounded uniformly for
$\delta$ in
$L^{\frac{\min\{\alpha+\theta_1,\gamma+\theta_2\}}{\alpha}}(Q_T)$
norm for any fixed $k>0$, we can use the Egrov theorem to conclude
that the last three terms on the right hand side of
(\ref{3-IV2+1}) vanish. Then we have
 \be\label{3-IV2+2}
\begin{split}
IV_2 \ge&\int_0^t\psi\int_\Omega\phi
\overline{T_k(\rho)}\,\,\overline{n^\alpha}\,dx\,ds.
\end{split}
\ee

{\noindent\bf Case 2.}  $\alpha\in[1,\infty),\,\,
\gamma\in(\frac{9}{5},\infty)$, and $\frac{1}{c_0}\rho_0\le n_0\le
c_0\rho_0$.

In this case, we have (\ref{2-r6}). Then repeating the
corresponding steps in Case 1, we get (\ref{3-IV2+2}).

\bigskip

For $IV_3$, we have \be\label{t1}
\begin{split}
IV_3=&\lim\limits_{\delta\rightarrow0^+}\int_0^t\psi\int_\Omega\phi
T_k(n_\delta)\rho_\delta^\gamma\,dx\,ds\\
=&\lim\limits_{\delta\rightarrow0^+}\int_0^t\psi\int_\Omega\phi
T_k\big(Ad_\delta\big)\rho_\delta^\gamma\,dx\,ds+\lim\limits_{\delta\rightarrow0^+}\int_0^t\psi\int_\Omega\phi
\Big(T_k(n_\delta)-T_k\big(Ad_\delta\big)\Big)\rho_\delta^\gamma\,dx\,ds\\
=&\lim\limits_{\delta\rightarrow0^+}\int_0^t\psi\int_\Omega\phi
T_k\big(Ad_\delta\big)B^\gamma
d_\delta^\gamma\,dx\,ds+\lim\limits_{\delta\rightarrow0^+}\int_0^t\psi\int_\Omega\phi
T_k\big(Ad_\delta\big)\Big(\rho_\delta^\gamma-B^\gamma
d_\delta^\gamma\Big)\,dx\,ds\\&+\lim\limits_{\delta\rightarrow0^+}\int_0^t\psi\int_\Omega\phi
\Big(T_k(n_\delta)-T_k\big(Ad_\delta\big)\Big)\rho_\delta^\gamma\,dx\,ds.
\end{split}
\ee Similar to the proof of (\ref{3-IV1}), we have
\be\label{t}\begin{split}
&\lim\limits_{\delta\rightarrow0^+}\int_0^t\psi\int_\Omega\phi
T_k\big(A d_\delta\big)B^\gamma d_\delta^\gamma\,dx\,ds\\ &\ge
\int_0^t\psi\int_\Omega\phi
\overline{T_k\big(A d\big)}B^\gamma\overline{d^\gamma}\,dx\,ds\\
&=\lim\limits_{\delta\rightarrow0^+}\int_0^t\psi\int_\Omega\phi
\big[T_k\big(Ad_\delta\big)-T_k\big(n_\delta\big)\big]B^\gamma\overline{d^\gamma}\,dx\,ds+\int_0^t\psi\int_\Omega\phi
\overline{T_k\big(n\big)}B^\gamma\overline{d^\gamma}\,dx\,ds\\
&=\lim\limits_{\delta\rightarrow0^+}\int_0^t\psi\int_\Omega\phi
\big[T_k\big(Ad_\delta\big)-T_k\big(n_\delta\big)\big]B^\gamma\overline{d^\gamma}\,dx\,ds+\int_0^t\psi\int_\Omega\phi
\overline{T_k\big(n\big)}\overline{\rho^\gamma}\,dx\,ds\\&+\lim\limits_{\delta\rightarrow0^+}\int_0^t\psi\int_\Omega\phi
\overline{T_k\big(n\big)}\Big(B^\gamma
d_\delta^\gamma-\rho_\delta^\gamma\Big)\,dx\,ds.
\end{split}\ee

For $\alpha,\gamma>\frac{9}{5}$ and
$\alpha\in(\gamma-\theta_1,\gamma+\theta_2)$, we have
$\gamma<\alpha+\theta_1$. In this case, $d_\delta^\gamma$ and
$\overline{d^\gamma}$ are bounded uniformly for $\delta$ in
$L^{\frac{\min\{\alpha+\theta_1,\gamma+\theta_2\}}{\gamma}}(Q_T)$
norm. For $\alpha\in[1,\infty),\,\,
\gamma\in(\frac{9}{5},\infty)$, and $\frac{1}{c_0}\rho_0\le n_0\le
c_0\rho_0$, we have (\ref{2-r6}) which implies that
$d_\delta^\gamma$ and $\overline{d^\gamma}$ are bounded uniformly
for $\delta$ in
$L^{\frac{\max\{\alpha+\theta_1,\gamma+\theta_2\}}{\gamma}}(Q_T)$
norm. Then using some similar arguments as in the estimates of
$IV_2$, we conclude that the last two terms on the right hand side
of (\ref{t1}) and the first and the third terms on the right hand
side of (\ref{t}) vanish. Thus \be\label{3-IV3} IV_3\ge
\int_0^t\psi\int_\Omega\phi
\overline{T_k\big(n\big)}\overline{\rho^\gamma}\,dx\,ds.\ee

(\ref{3-sum}) combined with the estimates of $IV_i$, i=1,2,3,4,
i.e., (\ref{3-IV1}), (\ref{3-IV4}), (\ref{3-IV2+2}), and
(\ref{3-IV3}), we have \be\label{3-sum+1}\begin{split}
&\int_0^t\psi\int_\Omega\phi \Big[\overline{T_k(\rho)}
+\overline{T_k(n)}\Big]\big(\overline{\rho^\gamma+n^\alpha}\big)\,dx\,ds\\
\le&\lim\limits_{\delta\rightarrow0^+}\int_0^t\psi\int_\Omega\phi
\big[T_k(\rho_\delta)+T_k(n_\delta)\big]\big(\rho_\delta^\gamma+n_\delta^\alpha\big)\,dx\,ds.
\end{split}
\ee where we have used
$$\overline{\rho^\gamma}+\overline{n^\alpha}=\overline{\rho^\gamma+n^\alpha}.$$

(\ref{3-sum+1}) implies (\ref{3-c1}). The proof of the lemma is
complete.

\endpf

Since $\psi$ and $\phi$ are arbitrary, we immediately get
\begin{corollary}\label{3-cor5.6}
Let $(\rho_\delta,n_\delta)$ be the solutions constructed in Proposition
\ref{2-le:aweak solution}, and $(\rho,n)$ be the limit, then \bex
 \Big[\overline{T_k(\rho)}
+\overline{T_k(n)}\Big]\big(\overline{\rho^\gamma+n^\alpha}\big)
 \le
\overline{\big[T_k(\rho)+T_k(n)\big]\big(\rho^\gamma+n^\alpha\big)}\eex
a.e. on $\Omega\times(0,T)$.
\end{corollary}

\begin{flushleft}
    \textbf{Step 2: Strong convergence of $\rho_{\delta}$ and $n_{\delta}$}
    \end{flushleft}

Here, we want to show the strong convergence of $\rho_{\delta}$ and $n_{\delta}$. This allows us to have Proposition \ref{key pro for delta}.
As in \cite{Feireisl}, we define \bex
L_k(z)=\left\{\begin{array}{l} z\log z,\quad 0\le z\le k, \\
[3mm] z\log k+z\displaystyle\int_k^z\frac{T_k(s)}{s^2}\,ds, \quad
z\ge k,
\end{array}
\right. \eex satisfying \bex\begin{split} L_k(z)=\beta_kz-2k\ for\
all\ z\ge 3k,
\end{split}
\eex   where
$$
\beta_k=\log
k+\displaystyle\int_k^{3k}\frac{T_k(s)}{s^2}\,ds+\frac{2}{3}.
$$
We denote $b_k(z):=L_k(z)-\beta_kz$ where $b^\prime_k(z)=0$ for
all large $z$, and \be\label{3-id} b^\prime_k(z)z-b_k(z)=T_k(z).
\ee

Note that $\rho_\delta, n_\delta \in L^{2}(Q_T)$, $\rho,n\in
L^{2}(Q_T)$, and $u_\delta, u\in L^2(0,T;H_0^1(\Omega))$. By Lemma
\ref{main lemma}, we conclude that $(n_\delta,u_\delta)$,
$(\rho_\delta,u_\delta)$, $(n,u)$ and $(\rho,u)$ are the
renormalized solutions of (\ref{a3-equation})$_i$ and
(\ref{3-equation})$_i$ for $i=1,2$, respectively. Thus we have
\bex\begin{cases}
[b_k(f_\delta)]_t+\mathrm{div}\big[b_k(f_\delta)u_\delta\big]+\big[b_k^\prime(f_\delta)f_\delta
-b_k(f_\delta)\big]\mathrm{div}u_\delta=0\quad \mathrm{in}\
\mathcal{D}^\prime(Q_T),\\[2mm]
[b_k(f)]_t+\mathrm{div}\big[b_k(f)u\big]+\big[b_k^\prime(f)f
-b_k(f)\big]\mathrm{div}u=0\quad \mathrm{in}\
\mathcal{D}^\prime(Q_T),
\end{cases}
\eex where $f_\delta=\rho_\delta,n_\delta$ and $f=\rho,n$.
 Thanks to (\ref{3-id}) and $b_k(z)=L_k(z)-\beta_kz$,
we arrive at \bex\begin{cases}
[L_k(\rho_\delta)+L_k(n_\delta)]_t+\mathrm{div}\big[\big(L_k(\rho_\delta)+L_k(n_\delta)\big)u_\delta\big]+\big[T_k(\rho_\delta)+T_k(n_\delta)\big]\mathrm{div}u_\delta=0\quad \mathrm{in}\ \mathcal{D}^\prime(Q_T),\\
[L_k(\rho)+L_k(n)]_t+\mathrm{div}\big[\big(L_k(\rho)+L_k(n)\big)u\big]+\big[T_k(\rho)+T_k(n)\big]\mathrm{div}u=0\quad
\mathrm{in}\ \mathcal{D}^\prime(Q_T).
\end{cases}
\eex  This gives \be\label{3-last}\begin{split}
&[L_k(\rho_\delta)-L_k(\rho)+L_k(n_\delta)-L_k(n)]_t+\mathrm{div}\big[\big(L_k(\rho_\delta)+L_k(n_\delta)\big)u_\delta
-\big(L_k(\rho)+L_k(n)\big)u\big]\\&+\big[T_k(\rho_\delta)+T_k(n_\delta)\big]\mathrm{div}u_\delta-\big[T_k(\rho)+T_k(n)\big]\mathrm{div}u=0.
\end{split}
\ee

Taking $\phi_j$ as the test function of (\ref{3-last}), and
letting $\delta\rightarrow\infty$, we have
\be\label{3-last1}\begin{split} &\int_\Omega
[\overline{L_k(\rho)}-L_k(\rho)+\overline{L_k(n)}-L_k(n)]\phi_j\,dx\\&-\lim\limits_{\delta\rightarrow0^+}
\int_0^t\int_\Omega\big[\big(L_k(\rho_\delta)+L_k(n_\delta)\big)u_\delta
-\big(L_k(\rho)+L_k(n)\big)u\big]\cdot\nabla\phi_j\,dx\,ds\\&+\lim\limits_{\delta\rightarrow0^+}\int_0^t\int_\Omega
\Big(\big[T_k(\rho_\delta)+T_k(n_\delta)\big]\mathrm{div}u_\delta-\big[T_k(\rho)+T_k(n)\big]\mathrm{div}u\Big)\phi_j\,dx\,ds=0,
\end{split}
\ee where \be\label{3-test2}\begin{split}&\phi_j\in
C_0^\infty(\Omega),\quad \phi_j(x)\equiv1\ \mathrm{for}\
\mathrm{any}\ x\in\big\{x\in\Omega \big|
\mathrm{dist(x,\partial\Omega)\ge\frac{1}{j}}\big\},\
0\le\phi_j\le1,\\& |\nabla\phi_j|\le c_0j,\ \phi_j\rightarrow
1\quad as\ m\rightarrow\infty
\end{split}
\ee for some positive $c_0$ independent of $j$.

Letting $j\rightarrow\infty$ in \eqref{3-last1}, we gain
\be\label{3-last2}\begin{split} &\int_\Omega
[\overline{L_k(\rho)}-L_k(\rho)+\overline{L_k(n)}-L_k(n)]\,dx\\=&-\lim\limits_{\delta\rightarrow0^+}\int_0^t\int_\Omega
\Big(\big[T_k(\rho_\delta)+T_k(n_\delta)\big]\mathrm{div}u_\delta-\big[T_k(\rho)+T_k(n)\big]\mathrm{div}u\Big)\,dx\,ds.
\end{split}
\ee 

In view of Lemma \ref{3-le:3.7}, we have
 \begin{equation}
\begin{split}
\label{last step-000}
&-\lim\limits_{\delta\rightarrow0^+}\int_0^t\int_\Omega
\big[T_k(\rho_\delta)+T_k(n_\delta)\big]\mathrm{div}u_\delta\,dx\,ds
\\=&-\frac{1}{2\mu+\lambda}\lim\limits_{\delta\rightarrow0^+}\int_0^t\int_\Omega
\big[T_k(\rho_\delta)+T_k(n_\delta)\big]\big[(2\mu+\lambda)\mathrm{div}u_\delta-\rho_\delta^\gamma-n_\delta^\alpha\big]\,dx\,ds
\\&-\frac{1}{2\mu+\lambda}\lim\limits_{\delta\rightarrow0^+}\int_0^t\int_\Omega
\big[T_k(\rho_\delta)+T_k(n_\delta)\big]\big[\rho_\delta^\gamma+n_\delta^\alpha\big]\,dx\,ds\\
=&-\frac{1}{2\mu+\lambda}\int_0^t\int_\Omega
\psi_j\phi_j\big[\overline{T_k(\rho)}+\overline{T_k(n)}\big]\big[(2\mu+\lambda)\mathrm{div}u-\overline{\rho^\gamma+n^\alpha}\big]\,dx\,ds
\\&-\frac{1}{2\mu+\lambda}\lim\limits_{\delta\rightarrow0^+}\int_0^t\int_\Omega
(1-\psi_j\phi_j)\big[T_k(\rho_\delta)+T_k(n_\delta)\big]\big[(2\mu+\lambda)\mathrm{div}u_\delta-\rho_\delta^\gamma-n_\delta^\alpha\big]\,dx\,ds
\\&-\frac{1}{2\mu+\lambda}\lim\limits_{\delta\rightarrow0^+}\int_0^t\int_\Omega
\big[T_k(\rho_\delta)+T_k(n_\delta)\big]\big[\rho_\delta^\gamma+n_\delta^\alpha\big]\,dx\,ds,
\end{split}
\end{equation}
 where $\psi_j$ and $\phi_j$ are given by (\ref{2-test1}) and
(\ref{2-test2}) respectively. Letting $j\rightarrow\infty$ in \eqref{last step-000}, we
have \be\label{3-last5}
\begin{split}
&-\lim\limits_{\delta\rightarrow0^+}\int_0^t\int_\Omega
\big[T_k(\rho_\delta)+T_k(n_\delta)\big]\mathrm{div}u_\delta\,dx\,ds
\\
=&-\frac{1}{2\mu+\lambda}\int_0^t\int_\Omega
\big[\overline{T_k(\rho)}+\overline{T_k(n)}\big]\big[(2\mu+\lambda)\mathrm{div}u-\overline{\rho^\gamma+n^\alpha}\big]\,dx\,ds
\\&-\frac{1}{2\mu+\lambda}\lim\limits_{\delta\rightarrow0^+}\int_0^t\int_\Omega
\big[T_k(\rho_\delta)+T_k(n_\delta)\big]\big[\rho_\delta^\gamma+n_\delta^\alpha\big]\,dx\,ds.
\end{split}
\ee

In view of (\ref{3-last2}) and (\ref{3-last5}), we have
\bex\begin{split} &\int_\Omega
[\overline{L_k(\rho)}-L_k(\rho)+\overline{L_k(n)}-L_k(n)]\,dx
 \\
=&\frac{1}{2\mu+\lambda}\int_0^t\int_\Omega
\big(\overline{T_k(\rho)}
+\overline{T_k(n)}\big)\big(\overline{\rho^\gamma}+\overline{n^\alpha}\big)\,dx\,ds
\\&-\frac{1}{2\mu+\lambda}\lim\limits_{\delta\rightarrow0^+}\int_0^t\int_\Omega
\big[T_k(\rho_\delta)+T_k(n_\delta)\big]\big(\rho_\delta^\gamma+n_\delta^\alpha\big)\,dx\,ds
\\&+\int_0^t\int_\Omega
[T_k(\rho)-\overline{T_k(\rho)}+T_k(n)-\overline{T_k(n)}]\mathrm{div}u\,dx\,ds,
\end{split}
\eex  with Corollary \ref{3-cor5.6}, which gives
\be\label{3-last4}\begin{split} \int_\Omega
[\overline{L_k(\rho)}-L_k(\rho)+\overline{L_k(n)}-L_k(n)]\,dx
\le\int_0^t\int_\Omega
[T_k(\rho)-\overline{T_k(\rho)}+T_k(n)-\overline{T_k(n)}]\mathrm{div}u\,dx\,ds.
\end{split}
\ee

Here we are able  to control the right-hand side of \eqref{3-last4} as in the following lemma.
\begin{lemma}
\label{right side lemma}
 \be\label{3-c2} \begin{split} \lim\limits_{k\rightarrow\infty}\int_0^t\int_\Omega
[T_k(\rho)-\overline{T_k(\rho)}+T_k(n)-\overline{T_k(n)}]\mathrm{div}u\,dx\,ds=0.
\end{split}
\ee
\end{lemma}
\pf
Recalling that $T(z)\le z$ for all $z$, we have
\bex\begin{split}\|T_k(\rho)-\overline{T_k(\rho)}\|_{L^2(Q_T)}\le&
\liminf\limits_{\delta\rightarrow0^+}\|T_k(\rho)-
T_k(\rho_\delta)\|_{L^{2}(Q_T)}\\
\le&C\liminf\limits_{\delta\rightarrow0^+}\|\rho+
\rho_\delta\|_{L^{\gamma+\theta_2}(Q_T)}\\ \le& C,\end{split} \eex
where we have used the H\"older inequality, $\gamma+\theta_2\ge2$,
(\ref{3-hierho}), (\ref{3-lim}), and (\ref{3-lim1}). With the help
of this estimate, (\ref{3-lim}), and (\ref{3-lim1}), one deduces
 \be\label{3-a5}
\begin{split} &\big|\int_0^t\int_\Omega
[T_k(\rho)-\overline{T_k(\rho)}]\mathrm{div}u\,dx\,ds\big|
\\ \le&\int_{Q_t\cap\{\rho\ge k\}}
|T_k(\rho)-\overline{T_k(\rho)}|\,|\mathrm{div}u|\,dx\,ds+\int_{Q_t\cap\{\rho\le
k\}} |T_k(\rho)-\overline{T_k(\rho)}|\,|\mathrm{div}u|\,dx\,ds\\
\le&\|T_k(\rho)-\overline{T_k(\rho)}\|_{L^2(Q_T)}\|\mathrm{div}u\|_{L^2\big(Q_t\cap\{\rho\ge
k\}\big)}+\|T_k(\rho)-\overline{T_k(\rho)}\|_{L^2\big(Q_t\cap\{\rho\le
k\}\big)}\|\mathrm{div}u\|_{L^2(Q_T)}\\
\le&C\|\mathrm{div}u\|_{L^2\big(Q_t\cap\{\rho\ge
k\}\big)}+C\|T_k(\rho)-\overline{T_k(\rho)}\|_{L^2\big(Q_t\cap\{\rho\le
k\}\big)}.
\end{split}
\ee

Note  that $T_k(z)=z$ if $z\le k$, we have
 \be\label{3-a6}
\begin{split}
\|T_k(\rho)-\overline{T_k(\rho)}\|_{L^2\big(Q_t\cap\{\rho\le
k\}\big)}=&\|\rho-\overline{T_k(\rho)}\|_{L^2\big(Q_t\cap\{\rho\le
k\}\big)}\\
\le&\liminf\limits_{\delta\rightarrow0^+}\|\rho_\delta-T_k(\rho_\delta)\|_{L^2(Q_T)}
\\
=&\liminf\limits_{\delta\rightarrow0^+}\|\rho_\delta-T_k(\rho_\delta)\|_{L^2(Q_T\cap\{\rho_\delta>k\})}
\\
\le&2\liminf\limits_{\delta\rightarrow0^+}\|\rho_\delta\|_{L^2(Q_T\cap\{\rho_\delta>k\})}
\\
\le&2k^{1-\frac{\gamma+\theta_2}{2}}\liminf\limits_{\delta\rightarrow0^+}\|\rho_\delta\|_{L^{\gamma+\theta_2}(Q_T)}^\frac{\gamma+\theta_2}{2}\rightarrow0
\end{split}
\ee as $k\rightarrow\infty$, due to (\ref{3-hierho}) and the
assumption $\gamma>\frac{9}{5}$ such that $\gamma+\theta_2>2$.

By (\ref{3-a5}) and (\ref{3-a6}), we conclude  \be\label{3-a7}
\lim\limits_{k\rightarrow\infty}\int_0^t\int_\Omega
[T_k(\rho)-\overline{T_k(\rho)}]\mathrm{div}u\,dx\,ds=0. \ee
Similarly, we have \be\label{3-a8}
\lim\limits_{k\rightarrow\infty}\int_0^t\int_\Omega
[T_k(n)-\overline{T_k(n)}]\mathrm{div}u\,dx\,ds=0. \ee With
(\ref{3-a7}) and (\ref{3-a8}),  (\ref{3-c2}) follows.
\endpf

Note that (\ref{3-last4}) and (\ref{3-c2}), we have \be\label{3-csum}\begin{split}
\limsup\limits_{k\rightarrow\infty}\int_\Omega
[\overline{L_k(\rho)}-L_k(\rho)+\overline{L_k(n)}-L_k(n)]\,dx \le
0.
\end{split}
\ee

By the definition of $L(\cdot)$, it is not difficult to justify
that \be\label{3-last6}\begin{cases}
\lim\limits_{k\rightarrow\infty}\Big[\|L_k(\rho)-\rho\log\rho\|_{L^1(\Omega)}+\|L_k(n)-n\log n\|_{L^1(\Omega)}\Big]=0,\\[2mm]
\lim\limits_{k\rightarrow\infty}\Big[\|\overline{L_k(\rho)}-\overline{\rho\log\rho}\|_{L^1(\Omega)}
+\|\overline{L_k(n)}-\overline{n\log n}\|_{L^1(\Omega)}\Big]=0.
\end{cases}
\ee

Since $\rho\log\rho\le \overline{\rho\log\rho}$ and $n\log
n\le\overline{n\log n}$ due to the convexity of $z\mapsto z\log
z$, we have \be\label{3-csum+1}\begin{split} 0\le\int_\Omega
[\overline{\rho\log\rho}-\rho\log\rho+\overline{n\log n}-n\log
n]\,dx \le 0,
\end{split}
\ee where we have used (\ref{3-csum}) and (\ref{3-last6}). Thus we obtain \bex \overline{\rho\log\rho}=\rho\log\rho \quad
\mathrm{and }\;\;\,\ \overline{n\log n}=n\log n. \eex
It allows us to have the
 strong
convergence of $\rho_\delta$ and $n_\delta$ in $L^\gamma(Q_T)$ and
in $L^\alpha(Q_T)$ respectively. Therefore we proved
(\ref{3-laststep}).

\endpf

With Proposition \ref{key pro for delta}, the proof of Theorem \ref{th:1.1} can be done.

\section*{Acknowledgements}

A. Vasseur's research was supported in part by NSF grant
DMS-1614918. A part of the work was done when H.Wen visited
Department of Mathematics at the University of Texas at Austin,
Texas, USA in September 2016. He would like to thank the
department for its hospitality. H.Wen's research was partially
supported by the National Natural Science Foundation of China
$\#$11722104, $\#$11671150 and supported by GDUPS (2016). C. Yu's
research was supported in part by Professor Caffarelli's NSF grant
DMS-1540162. A part of the work was done when C. Yu visited School
of Mathematics, South China University of Technology, Guangzhou,
China. He would like to thank the department for its hospitality.


\begin{thebibliography}{99}


\bibitem{BLS} J.-W. Barrett, Y. Lu, E. Suli,
Existence of large-data finite-energy global weak solutions to a compressible Oldroyd-B model,  arXiv:1608.04229


\bibitem{BBCMT} M. Baudin, C. Berthon, F. Coquel, R. Masson, Q. H. Tran, A relaxation method for two-fluid flow models with hydrodynamic closure law, Numer. Math. 99 (2005) 411-440.
\bibitem{BCT}M. Baudin, F. Coquel, Q.H. Tran, A semi-implicit relaxation scheme for modeling two-fluid flow in a pipeline, SIAM J. Sci. Comput. 27 (2005) 914-936.









\bibitem{Br} C. E. Brennen, Fundamentals of Multiphase Flow, Cambridge Univ. Press, 2005.



\bibitem{BJ}

D. Bresch, P.-E. Jabin,
Global Existence of Weak Solutions for Compresssible Navier--Stokes Equations: Thermodynamically unstable pressure and anisotropic viscous stress tensor
,   arXiv:1507.04629



\bibitem{CG}
J. A. Carrillo, T. Goudon,  Stability and asymptotic analysis of a fluid-particle interaction model. Comm. Partial Differ. Eqs. 31(7–9), 1349-1379 (2006)



\bibitem{DL1} R. J. DiPerna, P.-L. Lions,  Ordinary differential equations, transport theory and Sobolev spaces, Invent. Math. 98 (3) (1989) 511-547.

\bibitem{DL} R. J. DiPerna, P.-L. Lions, On the Cauchy problem for Boltzmann equations: global existence and weak stability. Ann. of Math. (2) 130 (1989), no. 2, 321-366.


\bibitem{E}S. Evje, Weak solution for a gas-iquid model relevant for describing gas-kick oil wells. SIAM J. Math. Anal. 43, 1887-1922 (2011).
\bibitem{EK}S. Evje, K.H. Karlsen, Global existence of weak solutions for a viscous two-fluid model. J. Differ. Equ. 245(9), 2660-2703 (2008).
\bibitem{EK2} S. Evje, K.H. Karlsen, Global weak solutions for a viscous liquid-gas model with singular pressure law. Commun. Pure Appl. Anal. 8, 1867-1894 (2009).


\bibitem{EWZ} S. Evje, H. Wen, C. Zhu, On global solutions to the viscous liquid-gas model with unconstrained transition to single-phase flow. Math. Models Methods Appl. Sci. 27 (2017), no. 2, 323-346.


\bibitem{Feiresil1}
E. Feireisl, On compactness of solutions to the compressible
isentropic Navier-Stokes equations when the density is not square
integrable, Comment. Math. Univ. Carolinae 42 (1) (2001), 83-98.

\bibitem{Feireisl}
E. Feireisl, A. Novotny, H. Petzeltova, On the existence of
globally defined weak solutions to the Navier-Stokes equations, J.
Math. Fluid Mech. 3(2001) 358-392.

\bibitem{Feireisl2} E. Feireisl,  Dynamics of Viscous
Compressible Fluids, Oxford University Press, Oxford, 2004.





\bibitem{Hoff87} D. Hoff, Global existence for 1D, compressible, isentropic Navier-Stokes equations with large initial data. Trans. Amer. Math. Soc. 303 (1987), no. 1, 169-181.


\bibitem{H95JDE} D. Hoff,  Global solutions of the Navier-Stokes equations for multidimensional,
compressible flow with discontinuous initial data.  J. Differential Equations 120
(1995), 215-254.



\bibitem{H95} D. Hoff, Strong convergence to global solutions for
multidimensional flows of compressible, viscous fluids with
polytropic equations of state and discontinuous initial data. Arch.
Rational Mech. Anal. \textbf{132} (1995), 1-14.




\bibitem{H97}D. Hoff, Discontinuous solutions of the Navier-Stokes equations for multidimensional
heat-conducting flow.  Arch. Rational Mech. Anal. 139 (1997), 303-354.


\bibitem{Hu}X. Hu, Hausdorff dimension of concentration for isentropic compressible Navier-Stokes equations,    arXiv:1606.06825
\bibitem{I} M. Ishii, Thermo-Fluid Dynamic Theory of two-fluid Flow, Eyrolles, Paris, 1975.

\bibitem{I2}
M. Ishii, One-Dimensional Drift-Flux Model and Constitutive Equations for Relative Motion between Phases in Various two-fluid Flow Regimes, Tech. report, Argonne National Laboratory Report ANL 77-47, 1977.


\bibitem{Jiang-Zhang}S. Jiang, P. Zhang, Global spherically symmetry solutions of
the compressible isentropic Navier-Stokes equations, Comm. Math.
Phys. 215(2001), 559-581.


\bibitem{KS} A.V. Kazhikhov, V.V. Shelukhin, Unique global solution with respect to time of initial-boundary value problems for one-dimensional equations of a viscous gas. J. Appl. Math. Mech. 41 (1977), no. 2, 273¨C282.; translated from Prikl. Mat. Meh. 41 (1977), no. 2, 282--291(Russian).


\bibitem{Lions}
P. L. Lions, Mathematical Topics in Fluid Mechanics, vol. II,
Compressible Models, Clarendon Press, Oxford, 1998.

\bibitem{Lions2}  P. L. Lions,  On some challenging problems in nonlinear partial differential equations. Mathematics:
frontiers and perspectives, 121-135, Amer. Math. Soc., Providence, RI, 2000.



\bibitem{MN79} A. Matsumura, T. Nishida,  The initial value problem for the equations of
motion of compressible viscous and heat-conductive fluids. Proc. Japan Acad.
Ser. A Math. Sci. 55 (1979), 337-342.



\bibitem{MN80}
A. Matsumura, T. Nishida, The initial value problem for the equations of
motion of viscous and heat-conductive gases. J. Math. Kyoto Univ. 20 (1980),
67-104.

\bibitem{MN83} A. Matsumura, T. Nishida,  Initial-boundary value problems for the equations of
motion of compressible viscous and heat-conductive fluids. Comm. Math. Phys., 1983, 89, 445-464



\bibitem{MMMNPZ} D. Maltese, M. Michálek, P. Mucha, A. Novotn\'y, M. Pokorn\'y, and E. Zatorska,\emph{ Existence of weak solutions for compressible Navier-Stokes equations with entropy transport. }J. Differential Equations 261 (2016), no. 8, 4448-4485.


\bibitem{MV2} A. Mellet, A. Vasseur, Existence and uniqueness of global strong solutions for one-dimensional compressible Navier-Stokes equations. SIAM J. Math. Anal. 39 (2007/08), no. 4, 1344-1365.



\bibitem{MV3} A. Mellet, A. Vasseur, Asymptotic analysis for a Vlasov-Fokker-Planck/compressible Navier-Stokes system of equations. Comm. Math. Phys. 281 (2008), no. 3, 573-596.

\bibitem{S}D. Serre,  Solutions faibles globales des ¨¦quations de Navier-Stokes pour un fluide compressible. C. R. Acad. Sci. Paris. I Math. 303 (1986), no. 13, 639-642.

\bibitem{6} J. Simon, Nonhomogeneous
viscous incompressible fluids: existence of velocity, density and
pressure, SIAM J. Math. Anal., 21(5) (1990), 1093-1117.







\bibitem{W}G. B. Wallis, One-Dimensional two-fluid Flow, McGraw-Hill, New York, 1979.





\bibitem{Yao-Zhang-Zhu} L. Yao, T. Zhang, C. Zhu, Existence and asymptotic behavior of global weak solutions to a 2D viscous liquid-gas
two-fluid flow model, SIAM J. Math. Anal. 42 (2010), 1874-1897.



\bibitem{Z}N. Zuber, On the dispersed two-phase flow in the laminar flow regime, Chem. Engrg. Sci., 19 (1964), pp. 897-917.











\end{thebibliography}


\end{document}